\PassOptionsToPackage{pdfpagelabels=false}{hyperref}

\documentclass[twoside,11pt,final]{entics} 

\usepackage{aliascnt}
\usepackage{amsmath, amssymb}
\usepackage{thmtools}
\usepackage{tikz-cd}
\usepackage{lmodern, bbold, stmaryrd}
\usepackage{hyperref}
\usepackage[capitalise]{cleveref}
\usepackage{enticsmacro}
\usepackage{graphicx}
\usepackage[all]{xy}
\usepackage[T1]{fontenc}  
\newcommand{\N}{\mathbb{N}}
\newcommand{\n}{\mathbb{n}}

\newcommand{\I}{\mathbb{I}}

\newcommand{\C}{\mathbf{C}}

\newcommand{\V}{\mathbf{V}}
\newcommand{\Set}{\mathbf{Set}}
\newcommand{\Cat}{\mathbf{Cat}}
\newcommand{\Enr}{\mathbf{Enr}}
\newcommand{\EnrCat}{\mathbf{EnrCat}}
\newcommand{\MonCat}{\mathbf{MonCat}}
\newcommand{\End}{\mathbf{Endo}}

\newcommand{\Alg}{\mathbf{Alg}}
\newcommand{\CoAlg}{\mathbf{CoAlg}}

\newcommand{\op}{\mathbf{^{op}}}

\newcommand{\ul}[1]{\underline{#1}}

\newcommand{\depth}{\operatorname{depth}}

\newcommand{\xto}[1]{\xrightarrow{#1}}
\renewcommand{\phi}{\varphi}
\renewcommand{\epsilon}{\varepsilon}
\newcommand{\id}{\operatorname{id}}
\newcommand{\const}{\operatorname{const}}

\newcommand{\ev}{\operatorname{ev}}

\newcommand{\Fr}{\operatorname{Fr}}
\newcommand{\Cof}{\operatorname{Cof}}

\newcommand{\m}{\mathfrak{m}}

\DeclareMathSymbol{\mathinvertedexclamationmark}{\mathord}{operators}{'074}
\DeclareMathSymbol{\mathexclamationmark}{\mathord}{operators}{'041}
\makeatletter
\newcommand{\raisedmathinvertedexclamationmark}{%
  \mathord{\mathpalette\raised@mathinvertedexclamationmark\relax}%
}
\newcommand{\raised@mathinvertedexclamationmark}[2]{%
  \raisebox{\depth}{$\m@th#1\mathinvertedexclamationmark$}%
}
\makeatother

\newcommand{\invexcl}{\raisedmathinvertedexclamationmark}
\newcommand{\fromI}[1]{\raisedmathinvertedexclamationmark_{#1}}

\newcommand{\dash}{\text{-}}

\def\conf{MFPS 2025} 	
\volume{5}			
\def\volu{5}			
\def\papno{16}			
\def\mydoi{16583}
\def\lastname{{\fontsize{10pt}{10pt}\selectfont Mulder, North, P{\'e}roux}} 
\def\runtitle{Functoriality of Enriched Data Types} 
\def\copynames{L. Mulder, P. North, M. P\'eroux}    
\begin{document}
\begin{frontmatter}
  \title{Functoriality of Enriched Data Types} 						
  \author{Lukas Mulder\thanksref{a}\thanksref{lukas.mulder@ru.nl}}	
   \author{Paige Randall North\thanksref{b}\thanksref{p.r.north@uu.nl}}		
   \author{Maximilien P\'eroux\thanksref{c}\thanksref{peroux@msu.edu}}	
   \address[a]{Institute for Computing and Inforation Science\\ Radboud University\\				
    Nijmegen, Netherlands}
   \thanks[lukas.mulder@ru.nl]{Email: \href{mailto:lukas.mulder@ru.nl} {\texttt{\normalshape
        lukas.mulder@ru.nl}}} 
  \address[b]{Department of Mathematics and Department of Information and Computing Sciences\\Utrecht University\\
    Utrecht, Netherlands} 
  \thanks[p.r.north@uu.nl]{Email:  \href{mailto:p.r.north@uu.nl} {\texttt{\normalshape
      p.r.north@uu.nl}}}
  \address[c]{Department of Mathematics\\Michigan State University\\
    East Lansing, Michigan, United States} 
  \thanks[peroux@msu.edu]{Email:  \href{mailto:peroux@msu.edu} {\texttt{\normalshape
    peroux@msu.edu}}}

\begin{abstract} 
  In previous work, categories of algebras of endofunctors were shown to be enriched in categories of coalgebras of the same endofunctor, and the extra structure of that enrichment was used to define a generalization of inductive data types. 
  These generalized inductive data types are parametrized by a coalgebra $C$, so we call them $C$-inductive data types; we call the morphisms induced by their universal property $C$-inductive functions.

  We extend that work by incorporating natural transformations into the theory: given a suitable natural transformation between endofunctors, we show that this induces enriched functors between their categories of algebras which preserve $C$-inductive data types and $C$-inductive functions.

  Such $C$-inductive data types are often finite versions of the corresponding inductive data type, and we show how our framework can extend classical initial algebra semantics to these types. 
  For instance, we show that our theory naturally produces partially inductive functions on lists, changes in list element types, and tree pruning functions.
\end{abstract}
\end{frontmatter}

\section{Introduction}\label{sec:introduction}

\subsection{Motivation}
Inductive types and inductively defined functions play a central role in functional programming languages. 
Their categorical semantics arise from interpreting them as initial algebras of an endofunctor. 
Giving categorical semantics to types, and more generally to programs, offers us rigorous tools to reason about programs.
Conversely, concepts emerging in category theory often find their way into functional programming languages.

In prior work \cite{north2023coinductive}, the categorical semantics of inductive data types was expanded by applying a construction from Sweedler theory \cite{chase1969hopf} to the category of algebras of an endofunctor. 
A key result showed that this category is enriched in the category of coalgebras of the same endofunctor, leading to the definition of $C$-inductive functions -- inductively defined functions parameterized by a coalgebra $C$.
By considering algebras that admit a unique $C$-inductive function to any other algebra, we obtained $C$-initial algebras, which generalize inductive data types, and call them $C$-inductive types.

In this paper, we extend this work by incorporating natural transformations into the theory.
Natural transformations between endofunctors provide the categorical account of how one inductively defines functions between inductive types -- they induce functions between the initial algebras.
Thus, here we show that they do the same thing for $C$-inductive functions.

We explain this concretely in the following three examples: first, an example of a $C$-inductive type; second, an example of a natural transformation inducing a function between inductive types; and third, an example of a natural transformation inducing a function between $C$-inductive types.

\begin{example}
    The set of finite binary trees with labels in a monoid $(M,\bullet, e)$, denoted $T_M$, is the initial algebra for the endofunctor $F(X) := 1 + M \times X \times X$.
	Being initial entails that there exists a unique algebra morphism $T_M \to A$ for all algebras $A$.
    Now, consider the set $S_n$ containing all possible binary tree shapes of depth at most $n \in \N$.
    The set $S_n$ naturally carries a coalgebra structure which unpacks a tree shape into the triple $(e,\ell,r)$, where $\ell$ and $r$ are the shapes of subtrees.
    An example of an $S_n$-inductive function is given by the function $S_n \times T_M \to A$  which takes a shape and a binary tree, prunes the tree according to the shape, and then applies the algebra morphism $T_M \to A$.
    Since we are always pruning trees to a maximum depth of of $n$, we can restrict ourselves to the algebra $T_{M,n}$ consisting of binary trees of depth at most $n$.
    There exists a unique $S_n$-inductive function $S_n \times T_{M,n} \to A$ for every algebra $A$, justifying calling $T_{M,n}$ a $S_n$-initial algebra, which corresponds to a $S_n$-inductive type.
\end{example}

Introducing the parameterization via a coalgebra facilitates a preprocessing step, allowing us to control up to which point input is being considered. 
The elements of the coalgebra serve as witnesses of this, ensuring that the input does not exceed a specified size, thereby enabling us to generalize inductive data types by imposing size restrictions.

Another fundamental notion in categorical semantics for functional programming is the natural transformation.
Defining a function between inductive data types
can be seen categorically as providing a natural transformation between the endofunctors defining the types.
Such an inductive function is given by constructing an algebra morphism from the initial algebra of an endofunctor $G$ to the initial algebra of another $F$. 
To do so, the initial $F$-algebra is endowed with a $G$-algebra structure by employing a natural transformation $\mu \colon G \to F$. 
The desired inductive function then arises as the unique algebra morphism from the initial $G$-algebra.
\begin{example}
    Continuing our running example, consider the endofunctor given by $G(X) := 1 + X$, whose initial algebra is the natural numbers $\N$. 
    We seek to define a function $\N \to T_M$ that maps an integer $n$ to the \emph{perfect} binary tree (i.e., every node has $2$ children) of depth $n$ with each node labeled by $e \in M$.
    This function is constructed using the natural transformation $\mu_X \colon 1 + X \to 1 + M \times X \times X$, defined component-wise as $\mu_X(x) = (e,x,x)$. 
    We show in Definition~\ref{def:pullback} that we can use $\mu_{T_M}$ to equip $T_M$ with the $G$-algebra structure $1 + T_M \xto{\mu_{T_M}} 1 + M \times T_M \times T_M \to T_M$.
	The unique $G$-algebra morphism $\N \to T_M$ then realizes the desired function.
\end{example}

We now unify the previous two examples by combining the theory of $C$-inductive types with natural transformations.
In this example we see how natural transformations allow us to carry over $C$-inductive functions from one functor to another, while retaining the underlying set of the coalgebra $C$.

\begin{example}
    To conclude our running example, consider the set $\n = \{0, \dots, n\}$, which can be given a $G$-coalgebra structure by sending $i \mapsto i-1$.
    We can define the $\n$-inductive function $\n \times \N \to A$, which takes a pair of numbers $(i,j)$ and applies the algebra morphism $\N \to A$ to $\min(i,j)$.
    Observe the coalgebra $\n$ serves a role similar to that of $S_n$ when pruning trees, as they both allow us to control up to what point the input is considered.
	
	Using $\mu \colon G \to F$, we equip the $G$-coalgebra $\n$ with an $F$-coalgebra structure $\n \to 1+\n \xto{\mu_{\n}} 1+ M \times \n \times \n$, which is called the \emph{pushforward} of $\n$ (cf. Definition~\ref{def:pushforward}).
	We note the resulting $F$-coalgebra is isomorphic to the $F$-coalgebra $S_{n,\text{perf}}$ of perfect binary tree shapes of depth at most $n$, and hence the underlying sets contain the same elements.
    In Definition~\ref{def:familyoffunctions} we see how the natural transformation $\mu$ acts on $\n$-inductive functions: for example, here it takes the function $\n \times \N \to \N$ to the function $S_{n,\text{perf}} \times T_M \to T_M$.
    Following the same reasoning as before, we can restrict to $T_{M,n}$ to obtain a function $S_{n,\text{perf}} \times T_{M,n} \to T_M$, and we can postcompose with $T_M \to A$ to obtain $\phi_A \colon S_{n,\text{perf}} \times T_{M,n} \to A$ for any $F$-algebra $A$.
    Intuitively, $\phi_A$ takes a depth and a binary tree, prunes the binary tree to that depth and then applies the algebra morphism $T_M \to A$ to the result.

	Using the natural transformation $\mu$, we can leverage $C$-inductive functions on $F$-algebras using the more insightful $G$-algebras and coalgebras.
	The key advantage is that $\mu$ preserves the underlying set of the coalgebra $C$ when acting on a $C$-inductive function as in Definition~\ref{def:familyoffunctions}.
	In particular, the coalgebra $S_{n,\text{perf}}$ has the same underlying set as $\n$, ensuring the coalgebraic control remains well-understood.
	
    We remark the use of the coalgebra $S_{n,\text{perf}}$ is essential, since simply defining an algebra morphism $T_{M,n} \to A$ does not work in general, as seen when taking $A = T_M$.
    Moreover, the obtained function $\phi_A$ is unique for every $A$, making $T_{M,n}$ an $S_{n,\text{perf}}$-initial algebra.
    Since $S_{n,\text{perf}}$ is obtained from $\n$ it is easy to handle -- every element of $S_{n,\text{perf}}$ corresponds to an element of $\n$.
    At the same time, using $\mu$ we do obtain the non-trivial $S_{n,\text{perf}}$-inductive data type $T_{M,n}$, and $S_{n,\text{perf}}$-inductive functions defined on it come with extra control given by the manageable $S_{n,\text{perf}}$.
\end{example}

Incorporating natural transformations into the theory of $C$-inductive types and functions gives us tools to construct these types and functions from simpler ones.
These types often impose size constraints, and we theorize this could be used to reason about program termination and memory usage though that will be left for future work.
The $C$-inductive functions give us more control when compared to regular inductive functions, and utilizing natural transformations to generate these functions allows us to reduce the complexity of the coalgebras involved.

In this paper, we develop a more robust theory for $C$-inductive types. 
Previous work has adopted an abstract perspective, focusing on enriched hom-objects defined via a universal property. 
Here, however, we focus on $C$-inductive functions, as they admit concrete definitions in all examples considered. 
This approach not only provides more constructive proofs but also enhances the feasibility of implementing the theory within a programming language, given the explicit characterization of $C$-inductive functions.

\subsection{Contributions}
In previous work, it was shown that for any well-behaved lax monoidal endofunctor, the category of algebras of that endofunctor is enriched in the category of its coalgebras \cite[Thm. 31]{north2023coinductive}. 
In this paper, we unify this result with monoidal natural transformations, which demonstrates this enrichment is functorial.

To formalize this, we use the fibered category $\EnrCat$ in Definition~\ref{def:EnrichedCat}, where objects are pairs $(\C, \V)$, with $\C$ as a $\V$-enriched category. 
By \cite[Thm. 31]{north2023coinductive}, the pair $(\Alg^F, \CoAlg^F)$ is an object of this category.

To establish functoriality, we take a lax monoidal natural transformation $\mu \colon F \to G$ between endofunctors on $\C$ and construct functors $(\mu_!,\mu_*) \colon (\Alg^F, \CoAlg^F) \to (\Alg^G, \CoAlg^G)$ in Theorem~\ref{thm:pushingforwardisfunctorial}, 
as well as $(\mu^*,\mu_{\invexcl}) \colon (\Alg^G, \CoAlg^G) \to (\Alg^F, \CoAlg^F)$ in Corollary~\ref{thm:pullingbackisfunctorial}. 
These functors are induced by transformations of measurings, which are equivalent by Theorem~\ref{thm:enrichedfunctor}.
Since $C$-initial algebras play a key role in our study, we show that the functor $(\mu_!,\mu_*)$ respects these structures in Theorem~\ref{lem:leftadjoijntpreservesCinitial}. 

Additionally, we provide worked examples in Section~\ref{sec:examples}, with a particular focus on expanding those introduced in the introduction in Section~\ref{sec:treepruning}.
We also demonstrate how a lax monoidal natural transformation induces a transformation of partially inductive functions on lists to partially inductive functions on natural numbers in Section~\ref{sec:pullingbacklists}.

\subsection{Related work}
Algebra and coalgebra are fundamental tools in computer science, particularly for giving semantics to inductive data types and dynamical systems. 
Coalgebras are widely used to model state-based systems and behaviors over time (an overview of which can be found in \cite{jacobs1997tutorial}), while initial algebras provide semantics for inductive data types such as natural numbers and trees. 
The theory behind inductive data types and their categorical semantics has been well developed, and has resulted in the theory of $W$-types \cite{martin2021intuitionistic,moerdijk2000wellfounded}.

There has also been considerable interest in generalizing these ideas in the spirit of Sweedler's theory \cite{chase1969hopf}, originally developed for $k$-algebras and $k$-coalgebras, where $k$ is some field. 
Analogues of Sweedler's theory have been developed for monoids \cite{fox1976universal,hyland2017hopf}, modules \cite{batchelor2000measuring}, dg-algebras \cite{anel2013sweedler}, $\infty$-algebras over $\infty$-operads \cite{peroux2022coalgebraic}, monads \cite{mcdermott2022sweedler}, and modules in double categories \cite{aravantinos2024enriched}.
Among these settings, algebras over an endofunctor are in a sense the most simple, as they avoid the coherence conditions required by more structured settings such as monads or operads.
Hence, here Sweedler's theory is laid out in the most transparent form.

In many concrete cases, $C$-inductive types are a given by imposing size constraints on inductive types, which is legitimated by the coalgebra $C$ dictating to which extent the induction proceeds.
This is related to techniques which guarantee the termination of recursive functions by indexing the number of recursive calls remaining before the function (forcefully) terminates \cite{appel2001indexed,mcbride2015turing,sozeau2020metacoq}.
For our $C$-inductive functions, elements of a coalgebra $C$ can also provide such indexing and we suspect $C$-inductive types can be used to reason about program termination in a similar fashion.

\section{Overview of prior work}\label{sec:prelims}
In this section we give a brief exposition of the material found in \cite{mulder2024measuring}.
We omit examples, as many are supplied in \cite[Sec. 5]{mulder2024measuring}.
The punchline of this section is that for any sufficiently well-behaved endofunctor $F: \C \to \C$ the category of $F$-algebras is enriched in the category of $F$-coalgebras.
For a background on monoidal categories and enrichment, see \cite{kelly1982enriched}.

We will be stating many results without proof, for which we defer the reader to \cite{mulder2024measuring}.
Throughout this section we will fix a monoidal category $(\C, \otimes, \I)$ and a lax monoidal endofunctor $(F, \nabla, \eta) :(\C, \otimes, \I) \to (\C, \otimes, \I)$.
Later in this section, we will place some extra conditions on $\C$ and $F$.
We will denote the category of $F$-algebras by $\Alg$ and the category of $F$-coalgebras by $\CoAlg$,
and we will often denote algebras by $(A,\alpha), (B,\beta) \in \Alg$ and coalgebras by $(C,\chi), (D,\delta) \in \CoAlg$.
For the sake of consistency, we will use terminology from previous work \cite{north2023coinductive,mulder2024measuring}, calling $C$-inductive functions \emph{measurings} and $C$-inductive types \emph{$C$-initial algebras}.

Let us start with what we will think of as the enriched morphisms.
There are many different ways to view these morphisms; as morphisms indexed by a coalgebra, as partial morphisms or as measurings.
\begin{definition}[\cite{north2023coinductive}, Def. 18]\label{def:measuring}
    Let $(A,\alpha), (B,\beta) \in \Alg$ and let $(C,\chi) \in \CoAlg$.
    We call a morphism $\phi \colon C \otimes A \to B$ in $\C$ a \emph{measuring from $A$ to $B$ by $C$} if it makes the diagram
    \[\begin{tikzcd}[row sep = tiny]
        & {F(C) \otimes F(A)} & {F(C \otimes A)} & {F(B)} \\
        {C \otimes F(A)} \\
        & {C \otimes A} && B
        \arrow["{\nabla_{C,A}}", from=1-2, to=1-3]
        \arrow["{F(\phi)}", from=1-3, to=1-4]
        \arrow["\beta"', from=1-4, to=3-4]
        \arrow["{\chi \otimes \id_{F(A)}}", from=2-1, to=1-2]
        \arrow["{\id_C \otimes \alpha}", from=2-1, to=3-2]
        \arrow["\phi", from=3-2, to=3-4]
    \end{tikzcd}\]
    commute.
    The set of all measurings from $A$ to $B$ by $C$ is denoted $\m_C(A,B)$.
\end{definition}
Precomposing a measuring with a coalgebra or algebra morphism or postcomposing with an algebra morphism will again result in a measuring.
Hence, we have a functor
$
\m \colon \CoAlg\op \times \Alg\op \times \Alg \to \Set.
$
If we place some conditions on $\C$ and $F \colon \C \to \C$, we will see this functor is representable in each of its three arguments.

\begin{remark}\label{rem:measuringbyunit}
    The monoidal unit $\I$ carries a coalgebra structure through $\eta \colon \I \to F(\I)$.
    A measuring from $A$ to $B$ by $\I$ and an algebra morphism $A \to B$ are equivalent.
\end{remark}

Given two measurings $\phi \colon C \otimes A \to A'$ and $\psi \colon D \otimes A' \to A''$ we would like to compose them to obtain a measuring from $A$ to $A''$.
This is a measuring by the coalgebra $D \otimes C$, which is given a coalgebra structure through the composition
$
D \otimes C \xto{\delta \otimes \chi} F(D) \otimes F(C) \xto{\nabla_{D,C}} F(D \otimes C).
$
The composition is then given by $\psi \circ (\id_D \otimes \phi)$, which makes the diagram in Definition~\ref{def:measuring} commute.
We will denote this composition of measurings by
$
\circ_\m \colon m_D(B,T) \times  m_C(A,B) \to  m_{D \otimes C}(A,T).
$
Next, we define $C$-initial objects, which give one motivation for our theory.
\begin{definition}[\cite{north2023coinductive}, Def.~35]
    Given a coalgebra $C \in \CoAlg$, we call an algebra $A \in \Alg$ a \emph{$C$-initial algebra} if for all $B \in \Alg$ there exists a unique measuring
    $
    \phi \colon C \otimes A \to B.
    $
\end{definition}

In order to develop the theory further, we assume $\C$ and $F: \C \to \C$ satisfy the following extra conditions.
From now on, we ask that $\C$ is a symmetric and closed monoidal category, and that it is \emph{locally presentable}. 
Furthermore we ask that $F$ is \emph{accessible}.
The motivation for these assumptions is that with them, the free and cofree functors $\Fr \colon \C \to \Alg$ and $\Cof \colon \C \to \CoAlg$ always exist.
These functors are used when constructing the representing objects of $\m$.
Moreover, as remarked earlier, the categories $\Alg$ and $\CoAlg$ are also locally presentable whenever $\C$ is, hence complete and cocomplete.
Another consequence is that $F$ is not only lax monoidal, but also \emph{lax closed}, respecting the closed structure of $\C$.

The punchline is that $\Alg$ is enriched in $\CoAlg$.
\begin{theorem}[\cite{north2023coinductive}, Thm.~31]\label{thm:algisenriched}
    The category $\Alg$ is enriched, copowered, and powered over the symmetric monoidal category $\CoAlg$ respectively via
    $$
    \Alg\op \times \Alg \xrightarrow{\ul{\Alg}(\_,\_)} \CoAlg, \quad \CoAlg \times \Alg \xrightarrow{\_\triangleright\_} \Alg, \quad \CoAlg\op \times \Alg \xrightarrow{[\_,\_]} \Alg  .
    $$
\end{theorem}
Summarizing the above, we have that for $A,B \in \Alg$ and $C \in \CoAlg$
$$
\m_C(A, B) \cong \CoAlg(C, \ul{\Alg}(A, B)) \cong \Alg(A,[C, B]) \cong \Alg(C \triangleright A, B),
$$
and hence $\m \colon \CoAlg\op \times \Alg\op \to \Alg$ is representable in each of its three arguments.
All these natural isomorphisms will be very useful down the road.

Using the representing objects, we have different ways to take an algebra and construct a coalgebra from it.
Two of these we give special attention, and are defined below.
\begin{definition}[\cite{north2023coinductive}, Def.~29]
    Let $I \in \Alg$ denote the initial algebra.
    We define the functor 
    $(\_)^* \colon \CoAlg\op \to \Alg$ by $C \mapsto [C,I]$
    and call $C^*$ the \emph{dual algebra of $C$}.
    We also define the functor
    $
        (\_)^\circ \colon \Alg\op \to \CoAlg$ by $
        A \mapsto \ul{\Alg}(A,I)
    $
    and call $A^\circ$ the \emph{dual coalgebra of $A$}.
\end{definition}

Using the isomorphisms regarding the representing objects of $\m_C(A,B)$, we can give some general results regarding preinitial algebras, their duals and $C$-initial algebras.
\begin{lemma}[\cite{north2023coinductive}, Lem.~37]\label{lem:dualissubterminal}
	Given algebras $P,B \in \Alg$ such that $P$ is preinitial, the coalgebra $\ul{\Alg}(P, B)$ is subterminal.
\end{lemma}
\begin{proof}
	An algebra $P$ is preinitial if and only if $\Alg(P,B)$ contains at most one element. Dually, a coalgebra $S$ is subterminal if $\CoAlg(C,S)$ contains at most one element.
	Since we have the isomorphism
	$
	\CoAlg(C, \ul{\Alg}(P, B)) \cong \Alg(P, [C,B])
	$
	and $P$ is preinitial, we conclude our result.
\end{proof}

\begin{corollary}
    Given a preinitial algebra $P \in \Alg$, its dual coalgebra $P^\circ$ is subterminal.
    \end{corollary}
\begin{proof}
    The dual coalgebra $P^\circ$ is defined as $\ul{\Alg}(P, I)$, where $I$ is the initial algebra.
	By Lemma~\ref{lem:dualissubterminal} we conclude our result.
\end{proof}
Next up is a powerful general result which allows us to generate a large number of $C$-initial algebras for specific coalgebras $C$.
\begin{proposition}\label{lem:preinitialisCinitial}
    A preinitial algebra $P \in \Alg$ is $P^\circ$-initial.
\end{proposition}
\begin{proof}
    The aim is to show $\m_{P^\circ}(P,B) \cong 1$ for all $B \in \Alg$.
    Our first remark is that there exists at most one measuring from $P$ to $B$ by $P^\circ$, since
    $
    \m_{P^\circ}(P,B) \cong \Alg(P, [P^\circ, B])
    $
    and $P$ is preinitial, hence has at most one morphism out of it.
    Second, since $P^\circ$ is subterminal by the previous lemma,  $\CoAlg(P^\circ,P^\circ) \cong 1$.
    Moreover, we have the identification
    $$
    1 \cong \CoAlg(P^\circ,P^\circ) \cong \Alg(P, (P^\circ)^*) = \Alg(P, [P^\circ, I]) \cong \m_{P^\circ}(P,I),
    $$
    where $I$ is the initial algebra.
    Postcomposing the unique measuring $\phi \in \m_{P^\circ}(P,I)$ with the morphism $\fromI{B} \colon I \to B$ yields a measuring
    $
    \fromI{B} \circ\phi \in \m_{P^\circ}(P,B).
    $
    We have shown there exists a unique measuring from $P$ to $B$ by $P^\circ$ for all $B \in \Alg$ and conclude $P$ is $P^\circ$ initial.
\end{proof}
As remarked before, the power of an initial algebra $I$ is that for any other algebra $B$ there exists a unique morphism $I \to B$.
For any preinitial algebra $P$ we know that there can be at most one algebra morphism $P \to B$, but there are no guarantees of it existing.
The above result tells us we can circumvent this disadvantage by not considering algebra morphisms $P \to B$, but instead measurings
$P^\circ \otimes P \to B$.
This gives the advantages of an initial algebra to a much broader class of algebras, namely all preinitial algebras.

\section{Enriched functors between categories of algebras}\label{sec:enrichedfunctors}
We build on \cite{north2023coinductive,mulder2024measuring}, of which a summary is provided in Section~\ref{sec:prelims}.
We write $\m_C(A,B)$ for the set of measurings from algebras $A$ to $B$ by a coalgebra $C$, which defines a functor $\m \colon \CoAlg\op \times \Alg\op \times \Alg \to \Set$.
Moreover, we will assume the category $\C$ and endofunctors $F,G$ on $\C$ satisfy the conditions of Theorem~\ref{thm:algisenriched}, the most notable of which is that $F$ and $G$ carry the structure of lax monoidal functors, denoted by $(\nabla^F, \eta^F)$, $(\nabla^G, \eta^G)$, respectively.
We will also assume all that natural transformations $\mu$ mentioned are lax monoidal.

The goal of this section is to show that a lax monoidal natural transformation $\mu: F \to G$ gives rise to functors $\mu_! \colon \Alg^F \to \Alg^G$ and $\mu^* :\Alg^F \to \Alg^G$ which respect the enrichment of algebras in coalgebras.
To this end, we first define the category $\EnrCat$ in Section~\ref{sec:transformingmeasurings}, which as objects has pairs $(\C, \V)$, where $\C$ is a $\V$-enriched category.
By Theorem~\ref{thm:algisenriched} we know $(\Alg^F, \CoAlg^F)$ and $(\Alg^G, \CoAlg^G)$ are elements of $\EnrCat$.
The main result of this section is that a coherent transformation of $F$-measurings into $G$-measurings is equivalent to a morphism $(\Alg^F, \CoAlg^F) \to (\Alg^G, \CoAlg^G)$ in $\EnrCat$.
With this idea in hand we show different ways in which a natural transformation results in transformations of measurings, and by extension morphisms in $\EnrCat$, in Section~\ref{sec:embeddingmeasurings} and Section~\ref{sec:pushingandpullingmeasurings}.
We wrap this all up by showing the functorial nature of these constructions, giving functors from the category of (lax monoidal) endofunctors on $\C$ to $\EnrCat$.

\subsection{Transforming measurings}\label{sec:transformingmeasurings}
We start off by giving a definition for the category of enriched categories, $\EnrCat$, which is defined using the Grothendieck construction.
We will denote the 2-category of monoidal categories, lax monoidal functors and lax monoidal natural transformations by $\MonCat$, and for any $\V \in \MonCat$ we will denote the 2-category of $\V$-enriched categories, $\V$-enriched functors and $\V$-enriched natural transformations by $\V\dash\Cat$.
Consider the 2-functor
\begin{align*}
	\Enr \colon \MonCat &\to \Cat \\
	\V &\mapsto \V\dash\Cat\\
	(\pi: \V \to \V') &\mapsto (\pi_* \colon \V\dash\Cat \to \V'\dash\Cat)\\
	(\nu: \pi \to \pi' )&\mapsto (\nu_* \colon \pi_* \to \pi'_*).
\end{align*}
Applying the Grothendieck construction to this functor yields the following category.
\begin{definition}\label{def:EnrichedCat}
	Let $\EnrCat$ be the category which as elements has pairs $(\C, \V)$, where $\V$ is a monoidal category and $\C$ is a $\V$-enriched category.
	A morphism $(\rho, \pi):(\C, \V) \to (\C',\V')$ is given by a pair of functors
	$
		\pi \colon \V \to \V' \in \MonCat,
		\rho \colon \pi_*(\C) \to \C' \in \V'\dash\Cat,
	$
	where $\pi_* \colon \V\dash\Cat \to \V'\dash\Cat$ is the change of base functor induced by the lax monoidal functor $\pi$ and $\rho$ is a $\V'$-enriched functor.

	Given pairs of functors $(\rho, \pi), (\rho', \pi') \in \EnrCat$, its composition
	$
	(\C,\V) \xrightarrow{(\rho, \pi)} (\C',\V') \xrightarrow{(\rho', \pi')} (\C'',\V'')
	$
	is given by $(\rho'', \pi'_*\circ \pi_*)$, where $\rho'' \colon  (\pi'_*\circ \pi_*)(\C) \to \C''$ is given by the composition $\rho'\circ \rho$ on objects.
	On morphisms $\rho''$ is given by the composition
	$
	\rho'_{\rho(A),\rho(B)} \circ \pi'(\rho_{A,B}) \colon (\pi' \circ \pi )_*(\C(A,B)) \to  \C''((\rho'\circ \rho)(A),(\rho'\circ \rho)(B)).
	$
	Identities are inherited from the underlying categories.
\end{definition}
We can think of $\EnrCat$ as the category containing all enriched categories, indexed by their enrichment.
Note that $\EnrCat$ is a 2-category, though we will not make use of this fact in this paper.

This section is motivated by us wishing to give a bijective correspondence between morphisms $(\rho, \pi) \colon (\Alg^F, \CoAlg^F) \to (\Alg^G, \CoAlg^G) \in \EnrCat$ and a natural transformations $\m^F_C(A,B) \to \m^G_{\pi(C)}(\rho(A),\rho(B))$ which \emph{respects composition of measurings}.
\begin{definition}\label{def:measuringrespectcomposition}
	Let $\rho \colon \Alg^F \to \Alg^G$ be a functor, $\pi \colon \CoAlg^F \to \CoAlg^G$ be a lax monoidal functor and let
	$
	\Phi_{A,B,C} \colon \m^F_C(A,B) \to \m^G_{\pi(C)}(\rho(A),\rho(B))
	$
	constitute a natural transformation from $\m^F$ to $\m^G \circ (\pi \times \rho \times \rho)$.
	We say $\Phi$ \emph{respects composition of measurings} if the diagram
	\[\begin{tikzcd}[row sep = small]
		{\m^F_D(B,T) \times \m^F_C(A,B)} && {\m^G_{\pi(D)}(\rho(B),\rho(T)) \times \m^G_{\pi(C)}(\rho(A),\rho(B))} \\
		&& {\m^G_{\pi(D)\otimes \pi(C)}(\rho(A),\rho(T))} \\
		{\m^F_{D\otimes C}(A,T)} && {\m^G_{\pi(D\otimes C)}(\rho(A),\rho(T))}
		\arrow["{\circ^F_\m}"', from=1-1, to=3-1]
		\arrow["{\Phi_{B,T,D} \times\Phi_{A,B,C}}", from=1-1, to=1-3]
		\arrow["{\Phi_{A,T,D\otimes C}}"', from=3-1, to=3-3]
		\arrow["{\circ^G_\m}", from=1-3, to=2-3]
		\arrow["{\m^G_{(\nabla^\pi)}(\rho(B),\rho(T))}"', from=3-3, to=2-3]
	\end{tikzcd}\]
	commutes.
\end{definition}
Since $\Phi$ in the above definition is a natural transformation, it automatically respects identities by Lemma~\ref{lem:enrichedfunctorrespectsidentities}, so we need not state this explicitly.
Since we think of measurings as our generalized algebra morphisms, this consistent way of transforming measurings should be enough to obtain an equivalent morphism in $\EnrCat$.

\begin{proposition}\label{lem:enrichedfunctor2}
	Let $\rho \colon \Alg^F \to \Alg^G$ be a functor and $\pi \colon \CoAlg^F \to \CoAlg^G$ be a lax monoidal functor.
    A natural transformation
	$
	\Phi_{A,B,C} \colon \m^F_C(A,B) \to \m^G_{\pi(C)}(\rho(A),\rho(B))
	$
	which respects composition induces a morphism
	$
	(\rho, \pi) \colon (\Alg^F, \CoAlg^F) \to (\Alg^G, \CoAlg^G)
	$
	in $\EnrCat$.
\end{proposition}
One might wonder if the converse is also true.
This is the case, showing there is an equivalence between natural transformations $\Phi$ which respect composition and morphisms in $\EnrCat$.
\begin{lemma}\label{lem:enrichedfunctor1}
	A morphism
	$
	(\rho, \pi) \colon (\Alg^F, \CoAlg^F) \to (\Alg^G, \CoAlg^G)
	$
	in $\EnrCat$ induces a natural transformation
	$
	\Phi_{A,B,C} \colon \m^F_C(A,B) \to \m^G_{\pi(C)}(\rho(A),\rho(B))
	$
	which respects composition.
\end{lemma}
Combing the above results, we obtain the following theorem.
\begin{theorem}\label{thm:enrichedfunctor}
	Let $\rho \colon \Alg^F \to \Alg^G$ be a functor, $\pi \colon \CoAlg^F \to \CoAlg^G$ be a lax monoidal functor.
	There exists a bijective correspondence between natural transformations
	$
	\Phi_{A,B,C} \colon \m^F_C(A,B) \to \m^G_{\pi(C)}(\rho(A),\rho(B))
	$
	which respects composition of measurings and morphisms
	$
	(\rho, \pi) \colon (\Alg^F, \CoAlg^F) \to (\Alg^G, \CoAlg^G)
	$
	in $\EnrCat$.
\end{theorem}
This theorem confirms our intuition, and in the next sections we will see how we can apply this theorem in different cases.
We conclude this section with a corollary which cements our intuition about composition of measurings.
Since we can compose morphisms in $\EnrCat$ and these are in bijective correspondence to natural transformations $\Phi$, it follows we can also compose these natural transformations.
\begin{corollary}
	If two natural transformations
	$\Phi_{A,B,C} \colon \m^F_C(A,B) \to \m^G_{\pi(C)}(\rho(A),\rho(B))$
	and
	$\Psi_{A,B,C} \colon \m^G_C(A,B) \to \m^H_{\pi'(C)}(\rho'(A),\rho'(B))$
	respect composition, then their composite
	$$
	(\Psi \circ \Phi)_{A,B,C} \colon \m^F_C(A,B) \to \m^G_{\pi(C)}(\rho(A),\rho(B)) \to \m^H_{\pi'(\pi (C))}(\rho'(\rho(A)),\rho'(\rho(B)))
	$$
	respects composition as well.
\end{corollary}

\subsection{Embedding measurings}\label{sec:embeddingmeasurings}
In this section we will see our first application of the theory presented in the previous section.
The guiding intuition is that given an section $\mu \colon F \rightarrowtail G$, this should result in an embedding $\Alg^F \rightarrowtail \Alg^G$ which also respects the enrichment.
This shows that measurings can be embedded analogously to algebras, while also serving as a stepping stone to the more involved constructions in Section~\ref{sec:pushingandpullingmeasurings}.

Before we get there, we first need to demonstrate what we can do given a natural transformation $\mu \colon F \to G$.
\begin{definition}\label{def:pullback}
    Given a natural transformation $\mu \colon F \to G$, the \emph{pullback functor} is defined as
    $
        \mu^* \colon \Alg^G \to \Alg^F,
        (\alpha \colon G(A) \to A) \mapsto (\alpha \circ \mu_A \colon F(A) \to G(A) \to A).
    $
\end{definition}
As expected, there is a dual to this definition for coalgebras.
\begin{definition}\label{def:pushforward}
    Given a natural transformation $\mu \colon F \to G$, the \emph{pushforward functor} is defined as
    $
	\mu_* \colon \CoAlg^F \to \CoAlg^G,
	(\chi \colon C \to F(C)) \mapsto (\mu_C \circ \chi :C \to F(C) \to G(C)).
    $
\end{definition}
We also have that this functor is \emph{strict}.
\begin{proposition}\label{lem:pushforwardisstrict}
    Given a lax monoidal natural transformation $\mu \colon F \to G$, the pushforward functor $\mu_* \colon \CoAlg^F \to \CoAlg^G$ is a strict monoidal functor 
    $
    \mu_* \colon (\CoAlg^F, \otimes, (\I, \eta_F)) \to (\CoAlg^G, \otimes, (\I, \eta_G)),
    $
\end{proposition}
\begin{proof}
    By $\mu$ being lax monoidal and the $\mu_*$ not changing that carrier of the algebra, all coherence maps are given identities.
	This implies $\mu_*$ is strict.
\end{proof}

Since we have already done a lot of work in the previous section relating transformations of measurings to morphisms in $\EnrCat$, we already have everything in place to state the key result of this section.
\begin{theorem}
    Let $(\nu \colon G \to F, \mu\colon F \to G)$ be a pair of lax monoidal natural transformations
    such that $\nabla^F_{C,A}$ coequalizes $\id_F \otimes \nu$ and $(\nu \otimes \nu) \circ (\mu \otimes \id_G)$.
	Then
    $
    (\nu^*, \mu_*) \colon (\Alg^G, \CoAlg^G) \to (\Alg^F, \CoAlg^F)
    $
    is a morphism in $\EnrCat$.
\end{theorem}

Returning to our intuition of a section $\mu \colon F \to G$ resulting in an enriched embedding of categories, we state the following result.
\begin{corollary}\label{thm:embeddingmeasurings}
	Given a pair of lax monoidal natural transformations $(\nu \colon G \to F, \mu\colon F \to G)$
	such that $\nu \circ \mu = \id_F$, then
    $
    (\nu^*, \mu_*) \colon (\Alg^G, \CoAlg^G) \to (\Alg^F, \CoAlg^F)
    $
    is a morphism in $\EnrCat$.
\end{corollary}
\begin{proof}
	Since $\nu \circ \mu = \id_F$, any morphism coequalizes $\id \otimes \nu_A$ and $(\nu_C \otimes \nu_A) \circ (\mu_C \otimes \id)$ since the morphisms coincide.
	By Corollary~\ref{thm:embeddingmeasurings}
	$
    (\nu^*, \mu_*) \colon (\Alg^G, \CoAlg^G) \to (\Alg^F, \CoAlg^F)
    $
    is a morphism in $\EnrCat$.
\end{proof}
As promised, we wrap everything up in a functor from the category of endofunctors to $\EnrCat$.
\begin{corollary}\label{lem:sectionsgiveenrichedmorphisms}
	Let $\C$ be a locally presentable, closed symmetric monoidal category and let $\End(\C)_{\text{retr.}}$ denote the category of accessible lax monoidal endofunctors on $\C$.
	Morphisms $G \to F$ in $\End(\C)_{\text{retr.}}$ are given by pairs $(\nu, \mu)$ of lax monoidal transformations such that $\nu \circ \mu = \id_F$.
	There exists a functor
	\begin{align*}
		\End(\C)_{\text{retr.}}&\longrightarrow \EnrCat \\
		F &\longmapsto (\Alg^F, \CoAlg^F)\\
		(\nu, \mu) &\longmapsto (\nu^*, \mu_*)
	\end{align*}
\end{corollary}

\subsection{Pushing forward and pulling back measurings}\label{sec:pushingandpullingmeasurings}
The previous section provides results which confirm our intuition.
However, we did require a lot, asking for a pair of lax monoidal natural transformations $F \xto{\mu} G \xto{\nu} F$ such that their composition is the identity on $F$.
Nicer would be if we could simply consider \textit{any} lax monoidal natural transformation $\mu \colon F \to G$, hopefully resulting in a morphism $(\Alg^F, \CoAlg^F) \to (\Alg^G, \CoAlg^G)$ in $\EnrCat$.

To establish this, we turn to the left adjoint of the pullback functor and the right adjoint of pushforward functor, denoted by $\mu_!$ and $\mu_{\invexcl}$ respectively.
In this section we aim to show that given a lax monoidal natural transformation $\mu \colon F \to G$, we obtain two morphisms
$
(\mu_!, \mu_*) \colon (\Alg^F, \CoAlg^F) \to (\Alg^G, \CoAlg^G)
$
and
$
(\mu^*, \mu_{\invexcl}) \colon (\Alg^G, \CoAlg^G) \to (\Alg^F, \CoAlg^F)
$
in $\EnrCat$.
\begin{remark}
One might wonder if instead of the left adjoint of $\mu^*$ we could also have considered the right adjoint.
Sadly the right adjoint does not exist, since in general $\mu^*$ does not preserve colimits.
Similarly, the left adjoint to $\mu_*$ does not exist since it does not preserve limits.
\end{remark}
Before we can get anything done, we must show the left adjoint to the pullback functor exists.
This could be proven using some adjoint functor theorem, but we prefer to give an explicit construction.
\begin{theorem}\label{thm:leftadjointofpullback}
    Given a natural transformation $\mu \colon F \to G$, the pullback functor $\mu^* \colon \Alg^G \to \Alg^F$ has a left adjoint $\mu_! \colon \Alg^F \to \Alg^G$ given by the coequalizer in $\Alg^G$,
	\[\begin{tikzcd}[ampersand replacement=\&]
		{\Fr^G(F(A))} \& {\Fr^G(A)} \& {\mu_!(A)},
		\arrow["{\tilde{f}}"', shift right=2, from=1-1, to=1-2]
		\arrow["{\Fr^G(\alpha)}", shift left=2, from=1-1, to=1-2]
		\arrow["{\text{coeq}}", dashed, from=1-2, to=1-3]
	\end{tikzcd}\]
	for any algebra $(A, \alpha) \in \Alg^F$.
    The morphism $\tilde{f}$ is obtained as adjunct under the free-forgetful adjunction of the composition
	$
	f \colon F(A) \xrightarrow{\mu_A} G(A) \xrightarrow{G(\eta_A)} G(\Fr^G(A)) \xrightarrow{\alpha_{\Fr^G}} \Fr^G(A)
	$
	with $\eta$ being the unit of the free-forgetful adjunction.
\end{theorem}
We would also like to exhibit the behavior of $\mu_!$ on morphisms.
Given an algebra morphism $g \colon A \to B \in \Alg^F$, we can draw the following diagram where we abbreviate $\Fr^G$ to $\Fr$:
\[\begin{tikzcd}
	{\Fr(F(A))} & {\Fr(A)} & {\mu_!(A)} \\
	{\Fr(F(B))} & {\Fr(B)} & {\mu_!(B)}.
	\arrow["{\Fr(F(g))}"', from=1-1, to=2-1]
	\arrow["{\Fr(\alpha)}", shift left=2, from=1-1, to=1-2]
	\arrow["{\Fr(\beta)}", shift left=2, from=2-1, to=2-2]
	\arrow["{\Fr(g)}", from=1-2, to=2-2]
	\arrow["{q_A}", from=1-2, to=1-3]
	\arrow["{q_B}"', from=2-2, to=2-3]
	\arrow["{\mu_!(g)}", dashed, from=1-3, to=2-3]
	\arrow["{\tilde{f}_B}"', shift right=2, from=2-1, to=2-2]
	\arrow["{\tilde{f}_A}"', shift right=2, from=1-1, to=1-2]
\end{tikzcd}\]
The morphism $\mu_!(g)$ is induced by the morphism $q_B \circ \Fr(g)$ coequalising $\Fr(\alpha)$ and $\tilde{f}_A$.

Next is a slightly surprising result which will turn out to be absolutely key throughout this entire section.
\begin{lemma}\label{lem:convolutionalgebracommutes}
    Let $F,G: \C \to \C$ be lax monoidal endofunctors and let $\mu \colon F \to G$ be a lax monoidal natural transformation.
	Let  $B\in \Alg^F$ and let $C \in \CoAlg^F$, then the $F$-algebras $[C,\mu^*(B)]$ and $\mu^*[\mu_*(C),B]$ are equal.
\end{lemma}
\begin{proof}
    Both $[C,\mu^*(B)]$ and $\mu^*[\mu_*(C),B]$ have $\underline{\C}(C,B)$ as underlying set.
	We claim $\id_{\underline{\C}(C,B)}$ is actually an algebra morphism.
	If this is the case, then $[C,\mu^*(B)]$ and $\mu^*[\mu_*(C),B]$ are equal.
	First, note $\mu$ is a lax monoidal natural transformation, hence a closed natural transformation \cite{eilenberg1966closedcat}.
	By definition of a closed natural transformation
	$
	\underline{\C}(F(C), \mu_B) \circ \tilde{\nabla}^F = \underline{\C}(\mu_C, G(B)) \circ \tilde{\nabla}^G \circ \mu_{\underline{\C}(C,B)}.
	$
	Verifying $\id_{\underline{\C}(C,B)}$ is an algebra morphism then amounts to a diagram chase.
\end{proof}
Using the Yoneda embedding and Lemma~\ref{lem:convolutionalgebracommutes}, we obtain the following result which will be used later to show $\mu_!$ preserves $C$-initial algebras.
\begin{lemma}\label{lem:leftadjointcommuteswithtensor}
	For all coalgebras $C \in \CoAlg^F$ and algebras $A \in \Alg^F$ there exists a natural isomorphism $\mu_*(C) \triangleright \mu_!(A) \cong \mu_!(C \triangleright A)$.
\end{lemma}

Now we have enough tools at our disposal to start building towards the main result of this section.
We aim to show a lax monoidal natural transformation $\mu \colon F \to G$ results in a morphism $(\mu_!, \mu_*) \colon (\Alg^F, \CoAlg^F) \to (\Alg^G, \CoAlg^G)$ in $\EnrCat$.
As a starting point, we will again make use of a transformation of measurings $\Phi$.
\begin{definition}\label{def:familyoffunctions}
	We define
	$
	\Phi_{A,B,C} \colon \m^F_C(A,B) \to \m^G_{\mu_*(C)}(\mu_!(A), \mu_!(B))
	$
	as the composite
	$$
	\m^F_C(A,B) \cong \Alg^F(C\triangleright A, B) \xto{\mu_!} 
	\Alg^G(\mu_!(C\triangleright A), \mu_!(B)) \cong 
	\Alg^G(\mu_*(C) \triangleright \mu_!(A), \mu_!(B)) \cong
	\m^G_{\mu_*(C)}(\mu_!(A), \mu_!(B)).
	$$
\end{definition}
Having this transformation of measurings is nearly enough, but as seen in the previous sections we also need to ask it respects composition to prove our main result.
\begin{lemma}\label{lem:familyoffunctionsrespectscomposition}
	The natural transformation defined in Definition~\ref{def:familyoffunctions} respects composition.
\end{lemma}

Summarizing the above, we state the following theorem.
\begin{theorem}
    Given a lax monoidal natural transformation $\mu \colon F \to G$, we obtain a morphism
    $
    (\mu_!, \mu_*) \colon (\Alg^F, \CoAlg^F) \to (\Alg^G, \CoAlg^G)
    $
    in $\EnrCat$.
\end{theorem}
\begin{proof}
	By Lemma~\ref{lem:familyoffunctionsrespectscomposition}, the natural transformation defined in Definition~\ref{def:familyoffunctions} using $\mu \colon F \to G$ respects composition.
	Using Theorem~\ref{thm:enrichedfunctor}, we conclude we obtain a morphism
	$
    (\mu_!, \mu_*) \colon (\Alg^F, \CoAlg^F) \to (\Alg^G, \CoAlg^G)
    $
    in $\EnrCat$.
\end{proof}
We wish to wrap this up in a concise statement, and do so by proving the following functor is well-defined.
The proof is identical to the proof of Corollary~\ref{lem:sectionsgiveenrichedmorphisms}
\begin{corollary}\label{thm:pushingforwardisfunctorial}
	Let $\C$ be a locally presentable, closed symmetric monoidal category and let $\End(\C)$ denote the category of accessible lax monoidal endofunctors on $\C$ and lax monoidal natural transformations.
	There exists a functor
	\begin{align*}
		\End(\C)&\longrightarrow \EnrCat \\
		F &\longmapsto (\Alg^F, \CoAlg^F)\\
		\mu &\longmapsto (\mu_!, \mu_*).
	\end{align*}
\end{corollary}

Since the enriched functor consists $(\mu_!, \mu_*)$ consists of left adjoints, one might wonder if it preserves $C$-initial algebras.
This is indeed the case, as shown by the following result.
\begin{theorem}\label{lem:leftadjoijntpreservesCinitial}
	Given a lax monoidal natural transformation $\mu \colon F \to G$, the enriched functor 
	$
	(\mu_!, \mu_*) \colon (\Alg^F, \CoAlg^F) \to (\Alg^G, \CoAlg^G)
	$ 
	preserves $C$-initial algebras.
\end{theorem}
\begin{proof}
	Let $A$ be $C$-initial. We claim $\mu_!(A)$ is $\mu_*(C)$-initial.
	Since $A$ is $C$-initial $C\triangleright A$ is the initial $F$-algebra.
	Since $\mu_!$ preserves initial objects by being a left adjoint we know $\mu_!(C\triangleright A)$ is the initial $G$-algebra.
	Finally, using Lemma~\ref{lem:leftadjointcommuteswithtensor}, we see $\mu_!(C\triangleright A) \cong \mu_*(C) \triangleright \mu_!(A)$ is the initial $G$-algebra, hence $\mu_!(A)$ is $\mu_*(C)$-initial.
\end{proof}

We now continue on to the dual story, where the right adjoint to the pullback functor plays a central role.
There are some subtle differences which we point out, and these subtle differences are also the reason we can not immediately dualize the statements above.
We start out by giving an explicit construction of the right adjoint $\mu_{\invexcl}$.
\begin{theorem}\label{thm:rightadjointexists}
    Given a natural transformation $\mu \colon F \to G$, the pushforward functor $\mu_* \colon \CoAlg^F \to \CoAlg^G$ has a right adjoint $\mu_{\invexcl} \colon \CoAlg^G \to \CoAlg^F$ given by the equalizer in $\CoAlg^F$
	\[\begin{tikzcd}
		{\mu_{\invexcl}(C)} & {\Cof^F(C)} & {\Cof^F(G(C))}
		\arrow["{\tilde{f}}"', shift right=2, from=1-2, to=1-3]
		\arrow["{\Cof(\chi)}", shift left=2, from=1-2, to=1-3]
		\arrow["{\operatorname{eq}}", dashed, from=1-1, to=1-2]
	\end{tikzcd}\]
	for any coalgebra $(C, \chi) \in \CoAlg^G$, where $\tilde{f}$ is obtained as adjunct under the free-forgetful adjunction of the composition
	$
	f \colon \Cof^F(C) \xrightarrow{\chi_{\Cof}} F(\Cof^F(C)) \xrightarrow{F(\epsilon)} F(C) \xrightarrow{\mu_C} G(C)
	$
	with $\epsilon$ being the counit of the free-forgetful adjunction.
\end{theorem}
\begin{proof}
	The proof is completely dual to that of Theorem~\ref{thm:leftadjointofpullback}
\end{proof}
The functor $\mu_{\invexcl}$ is lax monoidal if and only if $\mu_*$ is strong monoidal by \cite{kelly2006doctrinaladj}.
Since $\mu_*$ is a strict monoidal functor by Proposition~\ref{lem:pushforwardisstrict}, it is in particular a strong monoidal functor and hence $\mu_!$ is a lax monoidal functor.

At this point a subtle difference manifests itself.
Recall that when defining the natural transformation in Definition~\ref{def:familyoffunctions} which resulted in the morphism 
$
(\mu_!, \mu_*) \colon (\Alg^F, \CoAlg^F) \to (\Alg^G, \CoAlg^G)
$
in $\EnrCat$, the natural transformation was defined by applying $\mu_!$, combined with a bunch of natural isomorphisms.
Here, we apply the pullback functor $\mu^*$, as well as a natural transformation.
\begin{definition}\label{def:familyoffunctions2}
	We define
	$
	\Phi_{A,B,C} \colon \m^G_C(A,B) \to \m^F_{\mu_{\invexcl}(C)}(\mu^*(A), \mu^*(B))
	$
	as the composite
	\begin{multline*}
		\Phi_{A,B,C} \colon \m^G_C(A,B) \cong \Alg^G(A, [C,B]) \xto{\Alg^G(A, [\epsilon,B])} \Alg^G(A, [\mu_*\circ\mu_{\invexcl}(C),B]) \xto{\mu^*} \\
		\Alg^F(\mu^*(A), \mu^*([\mu_*\circ\mu_{\invexcl}(C),B])) = \Alg^F(\mu^*(A), ([\mu_{\invexcl}(C),\mu^*(B)])) \cong \m^F_{\mu_{\invexcl}(C)}(\mu^*(A), \mu^*(B)),
	\end{multline*}
	where $\epsilon$ is the counit of the adjunction $\mu_* \dashv \mu_{\invexcl}$.
\end{definition}
Defining $\Phi$ in this way ensures it is well-defined, but does make for an opaque definition.
Making it explicit at the level of the underlying category $\C$, $\Phi$ sends a $G$-measuring $f \colon C \otimes A \to B$ to the $F$-measuring 
$
f \circ (\epsilon_C \otimes \id_A) \colon \mu_{\invexcl}(C) \otimes A \to B,
$
where we ignore the pushforward functor $\mu_*$ and pullback functor $\mu^*$ since they do not change the carrier of the coalgebra or algebra.

Again, we need the defined natural transformation to respect composition.
\begin{lemma}\label{lem:familyoffunctionsrespectscomposition2}
	The natural transformation defined in Definition~\ref{def:familyoffunctions2} $\Phi$ respects composition.
\end{lemma}
Similar to before, we obtain the following result.
\begin{theorem}
    Given a lax monoidal natural transformation $\mu \colon F \to G$, we obtain a morphism
    $
    (\mu^*, \mu_{\invexcl}) \colon (\Alg^F, \CoAlg^F) \to (\Alg^G, \CoAlg^G)
    $
    in $\EnrCat$.
\end{theorem}
\begin{proof}
	By Lemma~\ref{lem:familyoffunctionsrespectscomposition2}, the natural transformation defined in Definition~\ref{def:familyoffunctions2} using $\mu \colon F \to G$ respects composition.
	Using Theorem~\ref{thm:enrichedfunctor}, we conclude we obtain a morphism
	$
    (\mu_!, \mu_*) \colon (\Alg^F, \CoAlg^F) \to (\Alg^G, \CoAlg^G)
    $
    in $\EnrCat$.
\end{proof}
We wrap up this section in the following statement.
Again, the proof is identical to the proof of Corollary~\ref{lem:sectionsgiveenrichedmorphisms}.
\begin{corollary}\label{thm:pullingbackisfunctorial}
	Let $\C$ be a locally presentable, closed symmetric monoidal category and let $\End(\C)$ denote the category of accessible lax monoidal endofunctors on $\C$ and lax monoidal natural transformations.
	There exists a functor
	\begin{align*}
		\End(\C)&\longrightarrow \EnrCat \\
		F &\longmapsto (\Alg^F, \CoAlg^F)\\
		\mu &\longmapsto (\mu^*, \mu_{\invexcl}).
	\end{align*}
\end{corollary}

\section{Examples}\label{sec:examples}
In this section we wish to exhibit the interplay between $C$-inductive functions and lax monoidal natural transformations.
Recall we call $C$-inductive functions \emph{measurings} in order to be consistent with \cite{north2023coinductive}.
Since measurings are our main object of interest, we will focus our attention on transformations of measurings which take a measuring for a functor $F$ and transforms it into a measuring for a functor $G$.
These transformations are equivalent to enriched functors constructed in the previous theorems by Theorem~\ref{thm:enrichedfunctor}.

Before turning to the examples, we provide an overview of the method used to construct the more involved among them.
There are two variants: one for pushing forward measurings and one for pulling them back.
First, given an endofunctor $F$, we identify the free algebra (respectively, cofree coalgebra) functor -- a step that often requires some creativity.
Using this (co)free functor, we then construct the left (respectively, right) adjoint to the pullback (respectively, pushforward) functor, via Theorem~\ref{thm:leftadjointofpullback} (respectively, Theorem~\ref{thm:rightadjointexists}).
Finally, we apply Lemma~\ref{lem:familyoffunctionsrespectscomposition} (respectively, Lemma~\ref{lem:familyoffunctionsrespectscomposition2}) to obtain the desired transformation of measurings.

We will often stress that a measuring is unique. 
Since the existence of a unique measuring implies any other measuring between the same objects coincides with this unique measuring, this allows us to easily prove two measurings coincide.
Morphisms out of the initial algebras share this uniqueness property, and was one of the motivations behind initial algebras as semantics for inductive data types.

Throughout this section we will use $(A,\alpha)$ and $(B,\beta)$ to denote algebras and $(C,\chi)$ to denote a coalgebra.
We may omit the algebra and coalgebra maps $\alpha, \beta, \chi$ if they are understood.
Moreover, we will denote the one-element set by $1 = \{*\}$.

\subsection{Natural numbers as lists}
As a first example we wish to exhibit the most simple transformation of measurings, given by an embedding of measurings as seen in Section~\ref{sec:embeddingmeasurings}.
The key observation is that natural numbers can be viewed as lists, where a number corresponds to the length of a list, regardless of its elements. 
We apply this perspective to arbitrary algebras and show that this embedding extends to measurings as well.

Let $(M, \bullet, e)$ be a monoid and consider the lax monoidal functors $F \colon X \mapsto 1 + X$ and $G \colon X \mapsto 1 + M \times X$.
The functor $F$ has $\N$ as initial algebra, and the functor $G$ has $M^*$ as initial algebra, where $M^*$ is the set of all finite lists with elements in $M$.
Define the lax monoidal natural transformations $\mu \colon F \to G$ by $\mu_X \colon 1 + X \to 1 + M \times X, x \mapsto (e,x)$ and $\nu \colon G \to F$ given by $\nu_X \colon 1 + M \times X \to 1 + X, (x', x) \mapsto x$.
We remark $\nu \circ \mu = \id_F$ in advance, allowing us to apply Corollary~\ref{thm:embeddingmeasurings} in the sequel.

We obtain the pullback functor $\nu^* \colon \Alg^F \to \Alg^G$, which maps an $F$-algebra $\alpha \colon 1 + A \to A$ to the $G$-algebra $1 + M \times A \xto{\mu} 1 + A \xto{\alpha} A$.
For example, equipping $F$-algebra $\N$ with the induced $G$-algebra structure $1 + M \times \N \to \N, * \mapsto 0, (x,i) \mapsto i+1$ captures the intuition that the specific elements of $M$ are irrelevant - only the fact that they contribute to list length matters.
Through $\nu^*$, the category $\Alg^F$ is embedded in $\Alg^G$ -- $\nu^*$ is a injective on objects and a fully faithful functor.

By Corollary~\ref{thm:embeddingmeasurings}, this embedding also carries over to measurings.
\begin{definition}
  The transformation of measurings $\Phi : \m^F_C(A,B) \to \m^G_{\mu_*(C)}(\nu^*(A), \nu^*(B))$ is given by sending a measuring $\phi \colon C \times A \to B$ to $\phi \colon \mu_*(C) \times \nu^*(A) \to \nu^*(B)$, which acts the same on the underlying objects $A$, $B$ and $C$.
\end{definition}
Intuitively, the conditions ensuring that $\phi$ is an $F$-measuring are sufficient for it to be a $G$-measuring when embedding $F$-algebras into $\Alg^G$ via $\nu^*$.

To establish that this embedding of measurings is well-defined, we provide the definitions of $F$-measurings and $G$-measurings, both special cases of Definition~\ref{def:measuring}. 
Given the similarity between $F$ and $G$, it is unsurprising their corresponding measurings closely resemble each other.
\begin{definition}
  Given $F$-algebras $\alpha \colon 1 + A \to A$ and $\beta \colon 1 + B \to B$ and an $F$-coalgebra $\chi \colon C \to 1 + C$, an $F$-measuring is a function $\phi \colon C \times A \to B$ such that
  \begin{enumerate}
      \item $\phi(c)(\alpha(*)) = \beta(*)$ for all $c \in C$
      \item $\phi(c)(\alpha(a)) = \beta(*)$ if $\chi(c) = *$
      \item $\phi(c)(\alpha(a)) = \beta(\phi(c')(a))$ if $\chi(c) = c'$.
  \end{enumerate}
\end{definition}
An $F$-measuring can be thought of as a partial inductively defined function, where the coalgebra $C$ dictates the extent to which induction proceeds.
\begin{definition}\label{def:listmeasuring}
  Given $G$-algebras $\alpha \colon 1 + M \times A \to A$ and $\beta \colon 1 + M \times B \to B$ and a $G$-coalgebra $\chi \colon C \to 1 + M \times C$, a $G$-measuring is a function $\phi \colon C \times A \to B$ such that
  \begin{enumerate}
      \item $\phi(c)(\alpha(*)) = \beta(*)$ for all $c \in C$
      \item $\phi(c)(\alpha(x,a)) = \beta(*)$ if $\chi(c) = *$
      \item $\phi(c)(\alpha(x,a)) = \beta(x' \bullet x, \phi(c')(a))$ if $\chi(c) = (x',c')$.
  \end{enumerate}
\end{definition}
Similarly, a $G$-measuring can be seen as a partial inductive function, though in this case the coalgebra $C$ also modifies the element $x \in M$ in condition (iii).
We check $\phi \colon \mu_*(C) \times \nu^*(A) \to \nu^*(B)$ is indeed a $G$-measuring by verifying the conditions given in Definition~\ref{def:listmeasuring}, using the fact that $\phi$ is an $F$-measuring.
\begin{enumerate}
  \item $\phi(c)(\alpha(\nu(*))) = \phi(c)(\alpha(*)) = \beta(*) = \nu(\beta(*))$
  \item $\phi(c)(\alpha(\nu(x,a))) = \phi(c)(\alpha(a)) = \beta(*) = \beta(\nu(*))$ if $\mu(\chi(c)) = \chi(c) = *$.
  \item $\phi(c)(\alpha(\nu(x,a))) = \phi(c)(\alpha(a)) = \beta(\phi(c')(a)) = \beta(\nu(e\bullet x, \phi(c')(a) ) )$ if $\mu(\chi(c)) = (e,c')$.
\end{enumerate}
Observing the above confirms the specific elements of $M$ are irrelevant to the embedded measuring, reinforcing our intuition regarding the embedded $F$-algebras.
This shows the embedding of $\Alg^F$ into $\Alg^G$ through $\nu^*$ respects measurings as well, and hence respects the enrichment of $\Alg^F$ and $\Alg^G$ in $\CoAlg^F$ and $\CoAlg^G$ respectively.

\subsection{Monoid homomorphisms inducing transformations of measurings}\label{ex:monoidhom}
In this example we consider the simplest non-trivial monoidal natural transformation, given by a monoid homomorphism.
This example also sheds light on the question posed in \cite{north2023coinductive}, namely how two different monoid structures on a set $M$ interact with each other in this context.

Let us set the scene first.
Take $(M, \bullet, e)$ to be a monoid and consider the lax monoidal functor $\const_M : \Set \to \Set, X \mapsto M$.
The algebras in question are functions $\alpha \colon M \to A$, and we can think of $\alpha$ as embedding elements of $M$ into $A$.
Coalgebras are given by $\chi \colon C \to M$, assigning to every element of $C$ and element of $M$.

We give a general definition of a measuring for any lax monoidal functor in Definition~\ref{def:measuring}, but would like to restate it here for our specific case.
\begin{definition}
    An $M$-measuring $\phi \colon C \times A \to B$ is a function $\phi$ such that $\phi(c, \alpha(x)) = \beta(\chi(c) \bullet x)$.
\end{definition}
We see a measuring takes an element $c\in C$ and an element $x \in M$ embedded in $A$, multiplies $\chi(c)$ and $x$ and embeds it in $B$.

Now that we have our bearings, we can work towards the transformation of measurings.
Given a monoid homomorphism $\mu \colon M \to M'$, we can define the lax monoidal natural transformation $\mu \colon \const_M \to \const_{M'}$.
In the sequel, the guiding intuition will be that we substitute elements of $M$ for elements of $M'$ using $\mu$.

By Theorem~\ref{thm:leftadjointofpullback} we can construct the left adjoint to the pullback functor $\mu^* \colon \Alg^{M'} \to \Alg^M$.
\begin{definition}
    The left adjoint to pullback functor, denoted $\mu_! : \Alg^M \to \Alg^{M'}$, sends an $M$-algebra $A$ to the pushout $A +_M M'$ of $\alpha \colon M \to A$ and $\mu \colon M \to M'$, given by
    $
    \mu_!(A) = A +_M M' = A + M' / \alpha(x) \sim \mu(x).
    $
    The algebra structure of $\mu_!(A)$ is given by $M' \ni x' \mapsto [x'] \in A +_M M'$.
\end{definition}
The functor $\mu_! : \Alg^M \to \Alg^{M'}$ takes an algebra $\alpha \colon M \to A$ and attempts to embed $M'$ into $A$ using the image of $\mu$.
Of course, if $\mu$ is not surjective this can not be done for all elements of $M'$, which is why the coproduct of $A$ and $M'$ is used.

We now have everything in place to utilize Lemma~\ref{lem:familyoffunctionsrespectscomposition} in order to obtain the following transformation of measurings.
\begin{definition}
    The transformation of measurings $\Phi : \m^M_C(A,B) \to \m^{M'}_{\mu_*(C)}(\mu_!(A), \mu_!(B))$ is given by sending a measuring $\phi \colon C \times A \to B$ to
    \begin{align*}
        \Phi(\phi) \colon C \times (A +_M M') &\to B +_M M'\\
        (c,[a]) &\mapsto [\phi(c,a)]\\
        (c,[x']) &\mapsto [\mu(\chi(c)) \cdot x'].
    \end{align*}
\end{definition}
Unpacking the above, the measuring $\Phi(\phi)$ acts the same as $\phi$ on $A$ in $A +_M M'$ and incorporates the new monoid $M'$ using $\mu$.
\begin{example}
    As a concrete example, consider the monoid homomorphism $\mu \colon (\{ \top, \bot \}, \wedge, \top) \to  (\{ \top, \bot \}, \vee, \bot)$ given by flipping truth values.
    In this case, an algebra $A$ of $(\{ \top, \bot \}, \wedge, \top)$ is given by assigning two truth values $\top_A, \bot_A \in A$.
    On the other hand, a coalgebra $C$ of $(\{ \top, \bot \}, \wedge, \top)$ assigns a truth value to every element of $C$.
    The functors $\mu_!$ and $\mu_*$ swap the interpretation of the two truth values of $A$ and the assigned truth values in $C$ respectively.
    The induced transformation of measurings follows the same pattern, again swapping the interpretation of the truth values.
\end{example}

\subsection{Pulling back lists}\label{sec:pullingbacklists}
In this example we expand on a standard way to define inductive functions by pulling back along natural transformations, as in Example~\ref{ex:pullbacklist}.
Again, let $(M, \bullet, e)$ be a monoid and consider the lax monoidal functors $F \colon X \mapsto 1 + X$ and $G \colon X \mapsto 1 + M \times X$.
Recall the functor $F$ has $\N$ as initial algebra, and the functor $G$ has $M^*$ as initial algebra, where $M^*$ is the set of all finite lists with elements in $M$.
Define the lax monoidal natural transformation $\mu \colon F \to G$ by $\mu_X \colon 1 + X \to 1 + M \times X, x \mapsto (e,x)$.
\begin{example}\label{ex:pullbacklist}
    As a starting point we wish to construct a function which takes a natural number $n$ and constructs a list of length $n$.
    Using $\mu$ we can pull back $M^*$ to the $F$-coalgebra $1 + M^* \xto{\mu_{M^*}} 1+ M \times M^* \to M^*$, giving us the unique $F$-algebra morphism $\N \to M^*, n \mapsto n *[e]$.
\end{example}

Since measurings are our main object of interest, we provide the definition for a $G$-measuring according to Definition~\ref{def:measuring}, as also stated in a previous example in Definition~\ref{def:listmeasuring}.
The definition of an $F$-measuring is given by replacing the arbitrary monoid $M$ for the trivial monoid $1$.
\begin{definition}
    Given $G$-algebras $\alpha \colon 1 + M \times A \to A$ and $\beta \colon 1 + M \times B \to B$ and a $G$-coalgebra $\chi \colon C \to 1 + M \times C$, a measuring is a function $\phi \colon C \times A \to B$ such that
    \begin{enumerate}
        \item $\phi(c)(\alpha(*)) = \beta(*)$ for all $c \in C$
        \item $\phi(c)(\alpha(x,a)) = \beta(*)$ if $\chi(c) = *$
        \item $\phi(c)(\alpha(x,a)) = \beta(x' \bullet x, \phi(c')(a))$ if $\chi(c) = (x',c')$.
    \end{enumerate}
\end{definition}
An important observation, and one of the most important features of a measurings, is that $\phi(c)(\alpha(x,a)) = \beta(*)$ if $\chi(c) = *$. 
In other words, if $c \in C$ has no successor, our inductive definition terminates.
This allows us to have control up to what point we are doing induction using elements of $C$.

In order to tie in the natural transformation $\mu$, we first need to compute the right adjoint of the pushforward functor $\mu_* : \CoAlg^F \to \CoAlg^G$ using Theorem~\ref{thm:rightadjointexists}.
\begin{definition}
    The right adjoint to the pushforward functor is given by $\mu_{\invexcl} : \CoAlg^G \to \CoAlg^F$, where
    $$
    \mu_{\invexcl}(C) = \{ c \in C \mid \chi(c) = * \text{ or } \chi(c) = (e,c'), c' \in \mu_{\invexcl}(C)\}.
    $$
\end{definition}
Note how $\mu_{\invexcl}(C) \subseteq C$, and consists of all elements which when unpacked give the neutral element or result in $*$.

We obtain the transformation of measurings which according to Lemma~\ref{lem:familyoffunctionsrespectscomposition2} is given by
\begin{align*}
    \Phi_{A,B,C} : \m^G_C(A,B) &\to \m^F_{\mu_{\invexcl}(C)}(\mu^*(A), \mu^*(B))\\
    \phi &\mapsto \phi | _{\mu_{\invexcl}(C) \times A}.
\end{align*}
We observe the transformation is given by restricting to $\mu_{\invexcl}(C)$, which intuitively carves out the part of a $G$-measuring which is also an $F$ measuring.

\begin{example}\label{ex:pullingbacklistmeasurings}
    As an example, consider the $G$-algebra $1 + M \times M^*_n \to M^*_n$ which as elements contains lists of length at most $n$.
    We also have the $G$-coalgebra ${M^*_n}^\circ \to 1 + M \times {M^*_n}^\circ$, where we use $^\circ$ to distinguish the coalgebra from the algebra.
    There exists a unique $G$-measuring $\phi \colon {M^*_n}^\circ \times M^*_n \to M^*$ which zips together two lists using the monoidal structure on $M$.
    The $G$-measuring $\phi$ gets sent to the $F$-measuring $\Phi(\phi) \colon {\{e\}^*_n}^\circ \times M^*_n \to M^*$, which zips a list containing only the element $e$ with a list of length at most $n$.
    Note we omit $\mu^*$ to ease on the notation.
    We could replace ${\{e\}_n^*}^\circ$ by $\n^\circ$, where $\n = \{0,\dots,n\}$, since we are always combining a elements of the list with the monoidal unit.
    From this new viewpoint, it becomes clear we are restricting a list in $M^*_n$ to a certain length $i \in \n$.
\end{example}

Tying this back to our initial idea of using a natural transformation to define a function $\N \to M^*$, we notice we can restrict ourselves to the $F$-algebra $\n \cong \{e\}^*_n \subseteq M^*_n$, similar to what we have done with the coalgebra.
This yields the unique measuring $\n^\circ \times \n \to M^*, (i,j) \mapsto \min(i,j) * [e]$.
It is important to stress the necessity of using a measuring here, since it is tempting to simply define a function $\n \to M^*, i \mapsto i * [e]$.
However, the defined function $\n \to M^*$ is not an $F$-algebra homomorphism.

One of the main motivations for using measurings is that we get guarantees about the size of the data structures involved.
Elements of the coalgebra $\n^\circ$ can be seen as a witness of the fact that this size can not exceed $n$, and hence are essential.
The idea to limit the length of the list by using $M^*_n$ instead of $M^*$ has the disadvantage of leaving control of the list to the algebra structure on $M^*_n$.
Instead, we use measurings to actively limit the length of the list when mapping into $M^*$, which gives us maximal control.

\subsection{Pruning trees}\label{sec:treepruning}

In this example we wish to exhibit to what degree measurings give us more flexibility and control than regular algebra homomorphisms.
We will showcase this by showing tree pruning naturally arises from measurings, and see that we even get more control over the resulting tree than simply pruning.
Introducing natural transformations will give us more tools to easily define measurings, and hence pruning functions.

Consider the monoid $(\N, + ,0)$ and the functors $G: X \mapsto 1 + \N \times X$ and $H:X \mapsto 1 + \N \times X \times X$.
The initial $G$-algebra is given by $\N^*$, the set of lists containing natural numbers, and the initial $H$-algebra is given by finite binary trees with nodes labeled by integers together with the empty tree, denoted $T = \{ (x,\ell, r) \mid x \in \N, \ell, r \in T \} \cup \{\emptyset\}$.
We will also be considering the final $H$-coalgebra, which is the set containing possibly infinite binary trees, denoted $T_\infty^\circ$, again using $^\circ$ to signal we are dealing with a coalgebra. 
Again, measurings are the central objects of study, so we apply the general definition of measuring Definition~\ref{def:measuring} to our specific case.
The definition of a $G$-measuring we have already seen in the previous example, and see the definition of an $H$-measuring follows a similar pattern.
\begin{defn}
    Given $H$-algebras $\alpha \colon 1 + M \times A \times A \to A$and $\beta \colon 1 + M \times B \times B \to B$ and a coalgebra $\chi \colon C \to 1 + M \times C \times C$, a measuring is given by a function $\phi \colon C \times A \to B$ such that
    \begin{enumerate}
        \item $\phi(c,\alpha(*)) = \beta(*)$ for all $c \in C$,
        \item $\phi(c, \alpha(x, a_\ell, a_r)) =  \beta(*)$ if $\chi(c) = *$,
        \item $\phi(c, \alpha(x, a_\ell, a_r)) = \beta(x' + x, \phi(c_\ell, a_\ell), \phi(c_\ell, a_\ell))$ if $\chi(c) = (x', c_\ell, c_e)$.
    \end{enumerate}
\end{defn}
We examine $H$-measurings in more detail to get a feel for them.
Intuitively, what a measuring does is take a tree-like $\alpha(x, a_\ell, a_r) \in A$ and overlaps it with the tree-like structure resulting from unfolding an element $c\in C$ using its coalgebra structure.
At all points where the tree nodes overlap we add the values in the nodes, and otherwise we discard the nodes.
Note there might be elements in $A$ which are not of the form $\alpha(x, a_\ell, a_r)$, and in that case we are free to choose where the measuring sends these elements.

\begin{definition}
    To ease the notation, we introduce the sum of trees with values in the natural numbers as
    \begin{align*}
        \oplus \colon T_\infty \times T &\to T \\
        (\emptyset, t) &\mapsto \emptyset\\
        (t, \emptyset) &\mapsto \emptyset\\
        ((x,\ell, r), (x',\ell', r')) &\mapsto (x + x', \ell \oplus \ell', r \oplus r')\\
    \end{align*}
\end{definition}
Note this is precisely the unique measuring from $T$ to $T$ by $T_\infty$, and many measurings considered will be restrictions or adaptations of this measuring.
This should come as no surprise, since the guiding intuition behind $H$-measurings is what this measurings is doing; overlapping trees and adding the values of the nodes.
\begin{example}
    If we restrict ourselves to the coalgebra $T^\circ_{\infty, \{0\}}$, possibly infinite binary trees containing only the value 0, we obtain the unique measuring $\oplus_0 \colon T^\circ_{\infty, \{0\}} \times T \to T$.
    In this case, the coalgebra $T^\circ_{\infty, \{0\}}$ can be thought of as the object containing all possible tree shapes.
    The unique measuring $\oplus_0$ is then given by $(t_0, t) \mapsto t_0 \oplus t$, hence prunes the tree $t$ according to the shape $t_0$.
    It does not alter the values contained in the tree $t$ since $t_0$ only contains the value $0$, the monoidal unit of $(\N, +, 0)$.

    We can extend this example to trees of finite depth.
    If we define $T_n$ to be the set of binary trees of depth at most $n$, and likewise for $T^\circ_n$ and $T_{n,\{0\}}^\circ$, we can consider the unique measuring $\oplus^n_0 \colon T_{n,\{0\}}^\circ \times T_n \to T$
    This gives us control over the depth of the trees, and hence the size of the trees we are considering.
\end{example}

Now that we have a feel for $H$-measurings, we can introduce a natural transformation and witness the interplay between $G$-measurings and $H$-measurings.
We define the lax monoidal natural transformation by
$
    \mu : G \to H ,
    \mu_X : 1 + \N \times X\ \to 1 + \N \times X \times X,
    * \mapsto *,
    (x,a) \mapsto (x,a,a).
$
In order to define a transformation from $G$-measurings to $H$-measurings, we need the left adjoint to the pullback functor $\mu^* \colon \Alg^H \to \Alg^G$, which can be constructed using Theorem~\ref{thm:leftadjointofpullback}.
\begin{definition}
    The left adjoint to pullback functor, $\mu_! : \Alg^G \to \Alg^H$, is given by
    $$
    \mu_!(A) = \left((M \times \mu_!(A) \times \mu_!(A)) + A \right) /\alpha(x,a) \sim (x,a,a).
    $$
    We denote elements as $[a], [x.\ell, r] \in \mu_!(A)$.
\end{definition}
Loosely speaking $\mu_!$ views list-like elements in $A$ as trees where all the nodes at the same level have the same value which we call \emph{equilevel trees}.
Following notation from Example~\ref{ex:pullingbacklistmeasurings}, let $\N^*_n$ denote the set of lists over $\N$ with length at most $n$.
Inspecting the definition of $\mu_!$, we notice the algebras $\N^*_n$ and $\N^*$ all get mapped to the initial algebra $T$ since their algebra morphisms are surjective.

Now we have everything in place to see how a transformation of measurings plays out.
By Lemma~\ref{lem:familyoffunctionsrespectscomposition}, we obtain the following transformation
\begin{definition}
    The of measurings $\Phi : \m^G_C(A,B) \to \m^H_{\mu_*(C)}(\mu_!(A), \mu_!(B))$ sends a $G$-measuring $\phi \colon C \times A \to B $ to the $H$-measuring
    \begin{align*}
        \Phi(\phi) \colon \mu_*(C) \times \mu_!(A) &\to \mu_!(B)\\
        \Phi(\phi)(c,[a]) &= [\phi(c,a)]\\
        \Phi(\phi)(c,[(x,\ell,r)]) &=
        \begin{cases}
            [\beta(*)] &\text{ if } \chi(c) = * \\
            [x'+ x, \Phi(\phi)(c',\ell),\Phi(\phi)(c',r)] &\text{ if } \chi(c) = (x',c').
        \end{cases}
    \end{align*}
\end{definition}
We see the transformation of measurings respects the intuition of list-like elements of $A$ and $B$ being equilevel trees.
It sends equilevel trees in $\mu_!(A)$ to the equilevel trees in $\mu_!(B)$ using $\phi$, and otherwise recurses on the separate branches until it finds an equilevel tree.
\begin{example}
    The transformation $\Phi$ sends the unique measuring $\phi \colon {\N^*_n}^\circ \times \N^*_n \to \N^*$ to $\Phi(\phi) \colon \mu_*({\N^*_n}^\circ) \times T \to T$, which prunes a tree according to the coalgebra $\mu_*({\N^*_n}^\circ)$.
    The coalgebra $\mu_*({\N^*_n}^\circ)$ contains lists of length at most $n$, but since we are pushing them forward by $\mu_*$ they should again be thought of as equilevel trees.
    In the context of the measuring $\Phi(\phi)$, the coalgebra $\mu_*({\N^*_n}^\circ)$ determines to what depth we are pruning the tree using the length of the list, but also gives us the option to alter values stored in the nodes.
    For instance, we can consider $\Phi(\phi)([0,1,2,\dots,i], t)$ which is the tree $t$ pruned to depth $i$, but where we also have added the level of the node to the value stored at each node.    
\end{example}

We wish to tie the above theory back to $C$-initial algebras and to do so focus our attention on the algebra $T_n$, which contain the trees of depth at most $n$.
\begin{lemma}
    The algebra $T_n$ is $\mu_*({\N^*_n}^\circ)$-initial.
\end{lemma}
\begin{proof}
    We can restrict the measuring $\Phi(\phi) \colon \mu_*({\N^*_n}^\circ) \times T \to T$ the algebra $T_n$.
    There exists at most one measuring out of $T_n$ to any target since $T_n$ is a preinitial algebra.
    We can explicitly construct this measuring using the fact that $T$ is the initial $H$-algebra and composition to obtain the measuring $\mu_*({\N^*_n}^\circ) \times T_n \to T \xto{\fromI{B}} B$.
    We conclude there exists a unique measuring $\mu_*({\N^*_n}^\circ) \times T_n \to B$ for every algebra $B$, hence that $T_n$ is $\mu_*({\N^*_n}^\circ)$-initial.
\end{proof}

As a final remark we would like to note we have restricted ourselves to the monoid $(\N,+,0)$. However, this idea works for any monoid $(M, \bullet, e)$.
One could even combine ideas from this example with that of Section~\ref{ex:monoidhom} and take some monoid homomorphism $h \colon M \to M'$ and construct the natural transformation
$ \mu_X \colon 1 + M \times X \to 1 + M' \times X \times X, * \mapsto *, (x,x') \mapsto (h(x), x', x')$.
This would still yield measurings which prune a tree with values in $M'$ up to a certain depth, and they also allow one to alter the values stored at the nodes of the tree using the using the monoid homomorphism $m$ and the monoid operation of $M'$.

\section{Conclusion \& Outlook}
We have shown the functorial nature of the enrichment of the category of algebras in the category of coalgebras.
In order to do so the category of enriched categories $\EnrCat$ in \ref{def:EnrichedCat}, of which $(\Alg^F, \CoAlg^F)$ are elements for a sufficiently well-behaved endofunctors by \cite[Thm. 31]{north2023coinductive}.
Functoriality was demonstrated in Corollary~\ref{thm:pushingforwardisfunctorial} and Corollary~\ref{thm:pullingbackisfunctorial} by constructing two functors
\begin{align*}
	\End(\C)&\longrightarrow \EnrCat \\
	F &\longmapsto (\Alg^F, \CoAlg^F)\\
	\mu &\longmapsto (\mu_!, \mu_*)
\end{align*}
and
\begin{align*}
	\End(\C)&\longrightarrow \EnrCat \\
	F &\longmapsto (\Alg^F, \CoAlg^F)\\
	\mu &\longmapsto (\mu^*, \mu_{\invexcl})
\end{align*}
where $\mu_!$ is the left adjoint of the pullback functor $\mu^*: \Alg^G \to \Alg^F$ and $\mu_{\invexcl}$ is the right adjoint of the pushforward functor $\mu_* \colon \CoAlg^F \to \CoAlg^G$.
We have taken measurings as our central object of study and provided explicit constructions of all functors involved, making the theory amenable to implementation in programming languages.

Since $C$-initial algebras exhibit categorical semantics similar to initial algebras, special significance has been assigned to them.
In Proposition~\ref{lem:preinitialisCinitial} we have shown any preinitial algebra $P$ is also $P^\circ$-initial, providing a broad class of $C$-initial algebras.
Furthermore, we proved in Theorem~\ref{lem:leftadjoijntpreservesCinitial} that the functor $(\mu_!, \mu_*)$ preserves $C$-initial algebras.

We have also seen the theory in action, most notably in an example regarding tree pruning.
In this example we have seen how measurings give us full control over the shape of the trees involved by pruning them, and natural transformations proved to be a useful tool to easily construct measurings.
More control over the shape and size of the data structure involved has been a general theme throughout our examples, and came naturally from using measurings.

Previous work introduced the notion of $n$-partial algebra homomorphisms \cite[Remark 4.7.5]{mulder2024measuring}. 
In future work, we aim to study these homomorphisms in greater depth and leverage them to construct a broader class of $C$-initial algebras. 
Additionally, we suspect that $(\mu_!, \mu_*)$ is left adjoint to $(\mu^*, \mu_{\invexcl})$ and plan to investigate this conjecture further. 
Another direction of interest is extending our framework beyond varying the endofunctor $F$, allowing for variations in the underlying category $\C$, which has so far remained fixed.
    
In this paper, we have taken a concrete approach using measurings, aiming to facilitate the implementation of this theory in a (possibly toy) programming language -- an avenue for future research. 
Having an implementation of measurings at our disposal should also enable the realization of $C$-initial algebras.
This introduces a new inductive principle guided by coalgebraic control, where the coalgebra dictates the extent of induction, akin to prior methods ensuring program termination \cite{appel2001indexed,mcbride2015turing,sozeau2020metacoq}.
By implementing $C$-initial algebras, we hope to further investigate the extent to which they can serve similar purposes.

\bibliographystyle{entics}
\bibliography{bibfile}

\appendix

\section{Details of \texorpdfstring{Section~\ref{sec:transformingmeasurings}}{Section 2.1}}\label{app:2.1}

\noindent
\textbf{Proposition 3.3}
\textit{	Let $\rho \colon \Alg^F \to \Alg^G$ be a functor and $\pi \colon \CoAlg^F \to \CoAlg^G$ be a lax monoidal functor.
    A natural transformation
	$
	\Phi_{A,B,C} \colon \m^F_C(A,B) \to \m^G_{\pi(C)}(\rho(A),\rho(B))
	$
	which respects composition induces a morphism
	$
	(\rho, \pi) \colon (\Alg^F, \CoAlg^F) \to (\Alg^G, \CoAlg^G)
	$
	in $\EnrCat$.}
	\smallskip

We will split the proof of Proposition~\ref{lem:enrichedfunctor2} over the following three lemmas.
First, we construct a morphism
$
\rho_{A,B} \colon \pi(\underline{\Alg}^F(A,B)) \to \underline{\Alg}^G(\rho(A),\rho(B))
$
in Lemma~\ref{lem:enrichedfunctorconstruction}, to obtain an enriched functor $\rho \colon \pi_*(\Alg^F) \to \Alg^G$ as in Definition~\ref{def:EnrichedCat}.
Second, we show $\rho_{A,B}$ respects identities in Lemma~\ref{lem:enrichedfunctorrespectsidentities}, and third we show $\rho_{A,B}$ respects composition in Lemma~\ref{lem:enrichedfunctorrespectscomposition}.

We start with the first step.
\begin{lemma}\label{lem:enrichedfunctorconstruction}
	Let $\rho \colon \Alg^F \to \Alg^G$ be a functor and $\pi \colon \CoAlg^F \to \CoAlg^G$ be a lax monoidal functor.
    A natural transformation
	$
	\Phi_{A,B,C} \colon \m^F_C(A,B) \to \m^G_{\pi(C)}(\rho(A),\rho(B))
	$
	induces a morphism
	$$
	\rho_{A,B} \colon \pi(\underline{\Alg}^F(A,B)) \to \underline{\Alg}^G(\rho(A),\rho(B)).
	$$
\end{lemma}
\begin{proof}
Consider the \emph{evaluation map}
$
\ev^F_{A,B} \colon \underline{\Alg}^F(A,B) \otimes A \to B,
$
which is the image of $\id_{\underline{\Alg}^F(A,B)}$ under the natural isomorphism
$
\m^F_{\underline{\Alg}^F(A,B)}(A,B) \cong
\CoAlg^F(\underline{\Alg}^F(A,B), \underline{\Alg}^F(A,B)).
$
Under $\Phi$ we obtain the $G$-measuring
$
\Phi(\ev^F_{A,B}) \colon \pi(\underline{\Alg}^F(A,B)) \otimes \rho(A) \to \rho(B).
$
Then ${\rho}_{A,B} \colon \pi(\underline{\Alg}^F(A,B)) \to \underline{\Alg}^G(\rho(A),\rho(B)) $ is obtained through the natural isomorphism
$
\m^G_{\pi(\underline{\Alg}^F(A,B))}(\rho(A), \rho(B)) \cong \CoAlg^G(\pi(\underline{\Alg}^F(A,B)), \underline{\Alg}^G(\rho(A),\rho(B))).
$
\end{proof}
Before we continue, we define the category of measurings from $A$ to $B$ to be the category which has measurings $\phi \colon C \otimes A \to B$ as objects and morphisms $f \colon \phi \to \phi'$ given by coalgebra morphisms $f \colon C \to C'$ such that $\phi = \phi' \circ (f \otimes \id_A)$.
Notice $\ev_{A,B} \colon \underline{\Alg}(A,B) \otimes A \to B$ is the terminal object in the category of measurings from $A$ to $B$ by the natural isomorphism $\m_C(A,B) \cong \CoAlg(C,\ul{\Alg}(A,B))$.
This will be the key to the following two proofs.

We now continue with the second step.
\begin{lemma}\label{lem:enrichedfunctorrespectsidentities}
	Let $\rho \colon \Alg^F \to \Alg^G$ be a functor, $\pi \colon \CoAlg^F \to \CoAlg^G$ be a lax monoidal functor and let $\Phi$ be
	a natural transformation
	$
	\Phi_{A,B,C} \colon \m^F_C(A,B) \to \m^G_{\pi(C)}(\rho(A),\rho(B)).
	$
	The induced morphism
	$
	\rho_{A,B} \colon \pi(\underline{\Alg}^F(A,B)) \to \underline{\Alg}^G(\rho(A),\rho(B))
	$
	from Lemma~\ref{lem:enrichedfunctorconstruction} respects identities.
\end{lemma}
\begin{proof}
We must show the following diagram commutes
\[\begin{tikzcd}[sep = tiny]
	&& {(\I, \eta_G)} \\
	& {\pi(\I, \eta_F)} \\
	{\pi(\underline{\Alg}^F(A,A))} &&&& {\underline{\Alg}^G(\rho(A),\rho(A))}
	\arrow["{\eta_\pi}"', from=1-3, to=2-2]
	\arrow["{{j^G_{\rho(A)}}}", from=1-3, to=3-5]
	\arrow["{{\pi(j_A^F)}}"', from=2-2, to=3-1]
	\arrow["{{{\rho}_{A,A}}}"', from=3-1, to=3-5]
\end{tikzcd}\]
where $j_A \colon \I \to \underline{\Alg}(A,A)$ is the family of identity elements of the enriched category.
We will turn to the category of measurings from $\rho(A)$ to $\rho(A)$ to show the above diagram commutes.
Writing $\lambda$ for the left unitor in the monoidal category $\C$, we have the measurings 
$((\I, \eta_G), \lambda_{\rho(A)})$
and 
$(\pi(\underline{\Alg}^F(A,A)), \Phi(\ev^F(A,A)))$.
The above diagram is a diagram of coalgebra morphisms, which precisely corresponds to a diagram of measurings from $\rho(A)$ to itself
\[\begin{tikzcd}[sep = small]
	& {((\I, \eta_G), \lambda_{\rho(A)})} \\
	{(\pi(\underline{\Alg}^F(A,A)), \Phi(\ev^F(A,A)))} && {(\underline{\Alg}^G(\rho(A),\rho(A)), \ev_{\rho(A),\rho(A)})}
	\arrow["{{\pi(j_A^F)} \circ \eta_\pi}"', from=1-2, to=2-1]
	\arrow["{{j^G_{\rho(A)}}}", from=1-2, to=2-3]
	\arrow["{{{\rho}_{A,A}}}"', from=2-1, to=2-3]
\end{tikzcd}\]
where ${\rho}_{A,A}$ and $j^G_{\rho(A)}$ are a morphism of measurings by definition, and $\pi(j_A^F) \circ \eta_\pi$ is morphism of measurings since $\pi$ is a lax monoidal functor, hence respects the left unitor.
The diagram commutes since $\underline{\Alg}^G(\rho(A),\rho(A))$ is the terminal object in the category of measurings from $\rho(A)$ to itself.
We conclude
$
j^G_{\rho(A)} = {\rho}_{A,A} \circ \pi(j_A^F)
$
and hence that $\rho$ respects identities.
\end{proof}
Notice that in the above two lemmas, we did not ask $\Phi$ to respect composition of measurings.
We conclude with the third step.
\begin{lemma}\label{lem:enrichedfunctorrespectscomposition}
	Let $\rho \colon \Alg^F \to \Alg^G$ be a functor, $\pi \colon \CoAlg^F \to \CoAlg^G$ be a lax monoidal functor and let $\Phi$ be
	a natural transformation
	$
	\Phi_{A,B,C} \colon \m^F_C(A,B) \to \m^G_{\pi(C)}(\rho(A),\rho(B))
	$
	which respects composition.
	The induced morphism
	$
	\rho_{A,B} \colon \pi(\underline{\Alg}^F(A,B)) \to \underline{\Alg}^G(\rho(A),\rho(B))
	$
	from Lemma~\ref{lem:enrichedfunctorconstruction} respects composition.
\end{lemma}
\begin{proof}
One can use a similar strategy as when showing $\rho_{A,B}$ respects identities for showing $\rho$ respects composition.
We must show the following diagram commutes
\[\begin{tikzcd}
	{\pi(\underline{\Alg}^F(B,T)) \otimes \pi(\underline{\Alg}^F(A,B))} & {\pi(\underline{\Alg}^F(A,T))} \\
	{\underline{\Alg}^G(\rho(B),\rho(T)) \otimes \underline{\Alg}^G(\rho(A),\rho(B))} & {\underline{\Alg}^G(\rho(A),\rho(T))}
	\arrow["{\pi(\circ^F)\circ \nabla^\pi}", from=1-1, to=1-2]
	\arrow["{\circ^G}", from=2-1, to=2-2]
	\arrow["{{\rho}_{B,T} \otimes {\rho}_{A,B}}"', from=1-1, to=2-1]
	\arrow["{{\rho}_{A,T}}"', from=1-2, to=2-2].
\end{tikzcd}\]
We turn to the category of measurings from $\rho(A)$ to $\rho(T)$, and recognize $(\underline{\Alg}^G(\rho(A),\rho(T)), \ev^G_{\rho(A),\rho(T)})$ is the terminal object in this category.
Moreover, we have the following measurings from $\rho(A)$ to $\rho(T)$
\begin{align*}
\big(\pi(\underline{\Alg}^F(B,T)) \otimes \pi(\underline{\Alg}^F(A,B)) ,\:& \Phi(\ev^F_{B,T}) \circ (\id \otimes \Phi(\ev^F_{A,B}))\big) \\
\big(\underline{\Alg}^G(\rho(B),\rho(T)) \otimes \underline{\Alg}^G(\rho(A),\rho(B)),\:& \ev^G_{\rho(B),\rho(T)} \circ (\id \otimes \ev^G_{\rho(A),\rho(B)}) \big)\\
\big(\pi(\underline{\Alg}^F(A,T)),\:& \Phi(\ev^F_{A,T})\big).
\end{align*}
We claim the above diagram of coalgebra morphisms corresponds exactly to a diagram of measurings in the category of measurings from $\rho(A)$ to $\rho(T)$.
The key step in verifying this is noting
$
\Phi(\ev^F_{B,T}) \circ (\id \otimes \Phi(\ev^F_{A,B})) = \Phi(\ev^F_{B,T} \circ (\id \otimes \ev^F_{A,B})) \circ \nabla^\pi,
$
since we asked $\Phi$ to respect composition of measurings.
Together with the naturality of $\Phi$ this implies $\pi(\circ^F)\circ \nabla^\pi$ is a morphism of measurings.
Since in the above diagram the composites map into the terminal object in the category of measurings from $\rho(A)$ to $\rho(T)$, they must coincide.
\end{proof}
Inspecting the proof above, we can state the following corollary which reduces the amount of verification needed to check if a natural transformation $\Phi$ respects composition.
\begin{corollary}\label{lem:respectseval}
	A natural transformation
	$ 
	\Phi_{A,B,C} \colon \m^F_C(A,B) \to \m^G_{\pi(C)}(\rho(A),\rho(B))
	$
	respects composition of measurings if and only if
	$
	\Phi(\ev^F_{B,T}) \circ (\id \otimes \Phi(\ev^F_{A,B})) = \Phi(\ev^F_{B,T} \circ (\id \otimes \ev^F_{A,B})) \circ \nabla^\pi.
	$
\end{corollary}

\noindent\textbf{Lemma 3.4}
\textit{	
	A morphism
	$
	(\rho, \pi) \colon (\Alg^F, \CoAlg^F) \to (\Alg^G, \CoAlg^G)
	$
	in $\EnrCat$ induces a natural transformation
	$
	\Phi_{A,B,C} \colon \m^F_C(A,B) \to \m^G_{\pi(C)}(\rho(A),\rho(B))
	$
	which respects composition.}
\smallskip

\begin{proof}
	We define the natural transformation $\Phi$ as
	\begin{multline*}
		\Phi \colon \m^F_C(A,B) \cong \CoAlg^F(C, \ul{\Alg^F}(A,B)) \xto{\pi} \\
		\CoAlg^G(\pi(C), \pi(\ul{\Alg^F}(A,B))) \xto{(\rho_{A,B})_*} \CoAlg^G(\pi(C), \ul{\Alg^G}(\rho(A),\rho(B))) \cong \m^G_{\pi(C)}(\rho(A),\rho(B)) .
	\end{multline*}
	By Corollary~\ref{lem:respectseval} it suffices to check
	$
	\Phi(\ev^F_{B,T}) \circ (\id \otimes \Phi(\ev^F_{A,B})) = \Phi(\ev^F_{B,T} \circ (\id \otimes \ev^F_{A,B})) \circ \nabla^\pi,
	$
	which is readily done via a diagram chase.
\end{proof}

\noindent\textbf{Theorem 3.5}
\textit{	
	Let $\rho \colon \Alg^F \to \Alg^G$ be a functor, $\pi \colon \CoAlg^F \to \CoAlg^G$ be a lax monoidal functor.
	There exists a bijective correspondence between natural transformations
	$
	\Phi_{A,B,C} \colon \m^F_C(A,B) \to \m^G_{\pi(C)}(\rho(A),\rho(B))
	$
	which respects composition of measurings and morphisms
	$
	(\rho, \pi) \colon (\Alg^F, \CoAlg^F) \to (\Alg^G, \CoAlg^G)
	$
	in $\EnrCat$.}
\smallskip

\begin{proof}
	Recall from the proof of Proposition~\ref{lem:enrichedfunctor2} the category of measurings from $A$ to $B$, which is defined as the category which has measurings $\phi \colon C \otimes A \to B$ as objects and morphisms $f \colon \phi \to \phi'$ given by coalgebra morphisms $f \colon C \to C'$ such that $\phi = \phi' \circ (f \otimes \id_A)$.
	The measuring $\ev_{A,B} \colon \underline{\Alg}(A,B) \otimes A \to B$ is the terminal object, and we denote the unique coalgebra morphism from a measuring $\phi$ to the terminal object by $!_\phi$.
	
	Recall that constructing the morphism $(\rho, \pi)$ from a natural transformation $\Phi$ is done by constructing the morphism $\rho_{A,B} \colon \pi(\ul{\Alg}^F(A,B)) \to \ul{\Alg}^G(\rho(A), \rho(B))$ using the terminal object $\ev_{\rho(A),\rho(B)}$.
	In the category of measurings from $\rho(A)$ to $\rho(B)$, we can draw the following commutative diagram.
	\[\begin{tikzcd}
		{\pi(C) \otimes\rho(A)} \\
		{\pi(\ul{\Alg}^F(A,B))} \otimes \rho(A) && {\rho(B)} \\
		{\ul{\Alg}^G(\rho(A), \rho(B))\otimes \rho(A)}
		\arrow["{\pi(!_\phi)\otimes \id_{\rho(A)}}"', from=1-1, to=2-1]
		\arrow["{\Phi(\phi)}", from=1-1, to=2-3]
		\arrow["{\Phi(\ev_{A,B})}"{description}, from=2-1, to=2-3]
		\arrow["{\rho_{A,B}\otimes \id_{\rho(A)}}"', from=2-1, to=3-1]
		\arrow["{\ev_{\rho(A), \rho(B)}}"', from=3-1, to=2-3]
	\end{tikzcd}\]
	The upper triangle commutes by $\Phi$ being a natural transformation, and the lower triangle commutes by the definition of $\rho_{A,B}$.
	Conversely, constructing a natural transformation
	$
	\Phi_{A,B,C} \colon \m^F_C(A,B) \to \m^G_{\pi(C)}(\rho(A),\rho(B))
	$
	from a morphism $(\rho, \pi)$ is done by sending a measuring $\phi \colon C\otimes A \to B$ to the bottom left composite in the diagram above.
	In this case, the diagram also commutes by definition.

	Since the diagram commutes in both cases, this shows there is a bijective correspondence between natural transformations $\Phi$ which respect composition of measurings and morphisms $(\rho, \pi)$ in $\EnrCat$.
\end{proof}

\section{Details of \texorpdfstring{Section~\ref{sec:embeddingmeasurings}}{Section 2.2}}\label{app:2.2}

\noindent\textbf{Theorem 3.10}
\textit{	
    Let $(\nu \colon G \to F, \mu\colon F \to G)$ be a pair of lax monoidal natural transformations
    such that $\nabla^F_{C,A}$ coequalizes $\id_F \otimes \nu$ and $(\nu \otimes \nu) \circ (\mu \otimes \id_G)$.
	Then
    $
    (\nu^*, \mu_*) \colon (\Alg^G, \CoAlg^G) \to (\Alg^F, \CoAlg^F)
    $
    is a morphism in $\EnrCat$.}
	\smallskip

\begin{proof}
    We define the natural transformation
    $
        \Phi_{A,B,C} \colon \m_C(A,B) \to \m_{\mu_*(C)}(\nu^*(A), \nu^*(B)),
        \phi \mapsto \phi.
    $
	This is well defined since given any $F$-measuring $\phi$, the following diagram
    \[\begin{tikzcd}[row sep = tiny, column sep = small]
        & {F(C) \otimes G(A)} & {G(C) \otimes G(A)} & {G(C \otimes A)} & {G(B)} \\
        {C \otimes G(A)} \\
        & {C \otimes F(A)} & {F(C) \otimes F(A)} & {F(C \otimes A)} & {F(B)} \\
        && {C\otimes A} && B
        \arrow["\phi"', from=4-3, to=4-5]
        \arrow["{\id \otimes\alpha}"', from=3-2, to=4-3]
        \arrow["\beta", from=3-5, to=4-5]
        \arrow["{\chi \otimes \id}", from=3-2, to=3-3]
        \arrow["{\nabla^F_{C,A}}", from=3-3, to=3-4]
        \arrow["{F(\phi)}", from=3-4, to=3-5]
        \arrow["{\id \otimes \nu_A}"', from=2-1, to=3-2]
        \arrow["{\chi \otimes \id}", from=2-1, to=1-2]
        \arrow["{\mu_C\otimes \id}", from=1-2, to=1-3]
        \arrow["{\nabla^G_{C,A}}", from=1-3, to=1-4]
        \arrow["{G(\phi)}", from=1-4, to=1-5]
        \arrow["{\nu_B}", from=1-5, to=3-5]
        \arrow["{\id \otimes \nu_A}"{description}, dashed, from=1-2, to=3-3]
        \arrow["{\nu_C \otimes \nu_A}"{description}, dashed, from=1-3, to=3-3]
        \arrow["{\nu_{C \otimes A}}"{description}, dashed, from=1-4, to=3-4]
    \end{tikzcd}\]
	commutes since $\nu$ is a lax monoidal natural transformation and $\nabla^F_{C,A}$ coequalizes $\id \otimes \nu_A$ and $(\nu_C \otimes \nu_A) \circ (\mu_C \otimes \id)$.
	This shows any $F$-measuring $\phi \colon C \otimes A \to B$ is also a $G$-measuring $\phi \colon \mu_*(C) \otimes \nu^*(A) \to \nu^*(B)$.
    Observe $\Phi$ respects composition by Corollary~\ref{lem:respectseval} and the fact that
    $$
	\Phi(\ev^F_{B,T}) \circ (\id \otimes \Phi(\ev^F_{A,B})) 
    = \ev^F_{B,T} \circ \id \otimes \ev^F_{A,B}
    = \Phi(\ev^F_{B,T} \circ (\id \otimes \ev^F_{A,B})) \circ \nabla^{\mu_*}
	$$
    since $\nabla^{\mu_*} = \id$ because $\mu_*$ is a strict monoidal natural transformation.
    By Theorem~\ref{thm:enrichedfunctor} we conclude 
    $
    (\nu^*, \mu_*) \colon (\Alg^G, \CoAlg^G) \to (\Alg^F, \CoAlg^F)
    $
    is a morphism in $\EnrCat$.
\end{proof}

\noindent\textbf{Corollary 3.12}
\textit{	
	Let $\C$ be a locally presentable, closed symmetric monoidal category and let $\End(\C)_{\text{retr.}}$ denote the category of accessible lax monoidal endofunctors on $\C$.
	Morphisms $G \to F$ in $\End(\C)_{\text{retr.}}$ are given by pairs $(\nu, \mu)$ of lax monoidal transformations such that $\nu \circ \mu = \id_F$.
	There exists a functor
	\begin{align*}
		\End(\C)_{\text{retr.}}&\longrightarrow \EnrCat \\
		F &\longmapsto (\Alg^F, \CoAlg^F)\\
		(\nu, \mu) &\longmapsto (\nu^*, \mu_*)
	\end{align*}}

\begin{proof}
	We have already defined the functor, and the only thing left to check is that it respects composition.
	Given two pairs of natural transformations $(\nu, \mu):H \to G$ and $(\nu', \mu') \colon G \to F$ , we obtain two morphisms $(\nu^*, \mu_*)$ and $(\nu'^*, \mu'_*)$.
	Composing the latter two morphisms in $\EnrCat$ yields the morphism $(\nu'^*, \mu'_*) \circ (\nu^*, \mu_*)$.
	By definition we already know $(\nu'^*, \mu'_*) \circ (\nu^*, \mu_*)$ and $((\nu' \circ \nu)^*,(\mu \circ \mu)_*)$ agree on objects.
	It remains to check the also agree on the enriched hom-objects.
	We need to verify the diagram
	\[\begin{tikzcd}[sep = tiny]
		{\mu'_*\circ \mu_* (\ul{\Alg^F}(A,B))} && {\ul{\Alg^H}(\nu'^*\circ \nu^* (A),\nu'^*\circ \nu^*(B))} \\
		& {\mu_*\ul{\Alg^G}(\nu^* (A), \nu^*(B))}
		\arrow["{(\nu' \circ \nu)^*_{A,B}}", from=1-1, to=1-3]
		\arrow["{\mu_*(\nu^*_{A,B})}"', from=1-1, to=2-2]
		\arrow["{{\nu'^*}_{\nu^*(A), \nu^*(B)}}"', from=2-2, to=1-3]
	\end{tikzcd}\]
	commutes.
	This is the case, since by definition all coalgebra morphisms involved are actually morphisms of measurings in the category of measurings from $\nu'^*\circ \nu^* (A)$ to $\nu'^*\circ \nu^*(B)$.
	Since $\ul{\Alg^H}(\nu'^*\circ \nu^* (A),\nu'^*\circ \nu^*(B))$ is the terminal object in this category, the morphisms must coincide.
	We conclude the functor respects composition.
\end{proof}

\section{Details of \texorpdfstring{Section~\ref{sec:pushingandpullingmeasurings}}{Section 2.3}}\label{app:2.3}

\noindent\textbf{Theorem 3.14} 
\textit{	
    Given a natural transformation $\mu \colon F \to G$, the pullback functor $\mu^* \colon \Alg^G \to \Alg^F$ has a left adjoint $\mu_! \colon \Alg^F \to \Alg^G$ given by the coequalizer in $\Alg^G$,
	\[\begin{tikzcd}[ampersand replacement=\&]
		{\Fr^G(F(A))} \& {\Fr^G(A)} \& {\mu_!(A)},
		\arrow["{\tilde{f}}"', shift right=2, from=1-1, to=1-2]
		\arrow["{\Fr^G(\alpha)}", shift left=2, from=1-1, to=1-2]
		\arrow["{\text{coeq}}", dashed, from=1-2, to=1-3]
	\end{tikzcd}\]
	for any algebra $(A, \alpha) \in \Alg^F$.
    The morphism $\tilde{f}$ is obtained as adjunct under the free-forgetful adjunction of the composition
	$
	f \colon F(A) \xrightarrow{\mu_A} G(A) \xrightarrow{G(\eta_A)} G(\Fr^G(A)) \xrightarrow{\alpha_{\Fr^G}} \Fr^G(A)
	$
	with $\eta$ being the unit of the free-forgetful adjunction.}
\smallskip

\begin{proof}
	Throughout this proof, we will write $\Fr^G = \Fr$ since we are only considering the free functor $\C \to \Alg^G$.
	Also note that we will omit the forgetful functor $U^G \colon \Alg^G \to \C$ to avoid a notational mess.
	
	Let $g \colon A \to B \in \C$. The idea is to show $g \in \Alg^F(A, \mu^*(B))$ if and only if its transpose $\tilde{g} \colon \Fr(A) \to B$ coequalizes $\Fr(\alpha)$ and $\tilde{f}$.
	If this is the case, then there is a one to one correspondence between algebra morphisms  $\Alg^F(A, \mu^*(B))$ and algebra morphisms $\Alg^G(\mu_{!}(A), B)$.

	First, assume $g \in \Alg^F(A, \mu^*(B))$. Since $\tilde{g}$ is the transpose of $g$ it is given by $\tilde{g} = \epsilon_B \circ \Fr(g)$.
	Similarly, $\tilde{f} = \epsilon_{\Fr(A)} \circ \Fr(f)$.
	Then show
	$
		\tilde{g} \circ\Fr( \alpha) = \epsilon_B \circ \Fr(\beta \circ G(g) \circ \mu_A)
	$
	using the previous equation and that $g$ is an algebra morphism.
	Next, use the triangle identities $\epsilon_{B} \circ \eta_{B} = \id_{B}$ and that $\tilde{g}$ is an algebra morphism to deduce
	$
	\tilde{g} \circ \alpha_{\Fr} = \beta \circ G(\tilde{g}).
	$
	Using this, continue to deduce
	$
		\tilde{g} \circ\Fr( \alpha) = \tilde{g}  \circ \tilde{f}\
	$
	using the triangle equalities, properties of natural transformations and transposes.
	Conclude $\tilde{g}$ coequalizes $\Fr(\alpha)$ and $\tilde{f}$ whenever $g \in \Alg^F(A,\mu^*(B))$.

	Second, assume $\tilde{g} \colon \Fr(A) \to B$ coequalizes $\Fr(\alpha)$ and $\tilde{f}$.
	Since $g \colon A \to B$ is the transpose of $\tilde{g}$, it is given by $g = \tilde{g} \circ \eta_A$.
	Using that $\tilde{g}$ coequalizes $\Fr(\alpha)$ and $\tilde{f}$ and properties of transposes we can deduce
	$
		g \circ \alpha = \tilde{g} \circ \epsilon_{\Fr(A)} \circ \eta_{\Fr(A)} \circ f.
	$
	Again, the triangle identities state $\epsilon_{\Fr(A)} \circ \eta_{\Fr(A)} = \id_{\Fr(A)}$.
	Using this, we can continue our deduction by
	$
		g \circ \alpha =\beta \circ \mu_B \circ F(g),
	$
	using the same properties as before.
	From this we deduce $g \in \Alg^F(A, \mu^*(B))$.

	We conclude the left adjoint of $\mu^*$ is given by the coequalizer of $\Fr(\alpha)$ and $\tilde{f}$.
\end{proof}

\noindent\textbf{Lemma 3.4}
\textit{	
	The natural transformation defined in Definition~\ref{def:familyoffunctions2} $\Phi$ respects composition.}

\begin{proof}
	We aim to use the fact that the Yoneda embedding is full and faithful in combination with Lemma~\ref{lem:convolutionalgebracommutes}.
	For any algebra $B \in \Alg^G$, we have the following natural isomorphisms
	\begin{align*}
		\Alg^G(\mu_*(C)\triangleright \mu_!(A), B)
		&\cong\Alg^G(\mu_!(A), [\mu_*(C), B])\\
		&\cong \Alg^F(A, \mu^*([\mu_*(C), B]))\\
		&= \Alg^F(A, [C, \mu_*(B)])\\
		&\cong \Alg^F(C\triangleright A, \mu^*(B))\\
		&\cong \Alg^G(\mu_!(C\triangleright A), B).
	\end{align*}
	Since the Yoneda embedding is full and faithful, we can conclude there is a natural isomorphism $\mu_*(C) \triangleright \mu_!(A) \cong \mu_!(C \triangleright A)$.
\end{proof}

\begin{lemma}\label{lem:tensorrespectsclosed}
	For all coalgebras $C, D \in \CoAlg$ and algebras $B \in \Alg$ we have the natural isomorphism
	$
	[C, [D,B]] \cong [D \otimes C , B ].
	$
\end{lemma}

\begin{proof}
	The underlying objects of these algebras are given by $\ul{\C}(C, \ul{\C} (D,B))$ and $\ul{\C}(D \otimes C, B)$ respectively.
	Since $\C$ is a closed monoidal category, we know there is a natural isomorphism
	$
	\ul{\C}(C, \ul{\C} (D,B)) \cong \ul{\C}(D \otimes C, B).
	$
	We claim this isomorphism lifts to an isomorphism of algebras and this can be verified by a diagram chase using the closed monoidal structure on $\C$ and naturality of the isomorphism $\ul{\C}(C, \ul{\C} (D,B)) \cong \ul{\C}(D \otimes C, B)$.
\end{proof}
	
\begin{lemma}\label{lem:tensorcommutes}
	For all coalgebras $C, D \in \CoAlg$ and algebras $A \in \Alg$ we have the natural isomorphism
	$$
	D \triangleright (C\triangleright A) \cong (D \otimes C) \triangleright A.
	$$
\end{lemma}

\begin{proof}
	We again aim to use the fact the the Yoneda embedding is full and faithful in combination with Lemma~\ref{lem:tensorrespectsclosed}.
	For any algebra $B \in\Alg$ we have the following natural isomorphisms
	\begin{align*}
		\Alg(D \triangleright (C\triangleright A) , B)
		&\cong \Alg(C \triangleright A, [D,B])\\
		&\cong \Alg(A, [C,[D,B]])\\
		&\cong \Alg(A, [C\otimes D,B])\\
		&\cong \Alg( (D \otimes C) \triangleright A, B)
	\end{align*}
	and by the Yoneda embedding we conclude our result.
\end{proof}

\noindent\textbf{Lemma 3.18}
\textit{	
	The natural transformation defined in Definition~\ref{def:familyoffunctions} respects composition.}

\begin{proof}
	We aim to show the diagram from Definition~\ref{def:measuringrespectcomposition} commutes.
	To do so, we remark that composition of measurings
	$
		\circ_\m :\m_D(B,T) \times \m_C(A,B) \to \m_{D\otimes C}(A,T),
		(\psi, f) \mapsto \psi \circ (\id_D \otimes f)
	$
	under the natural identifications $\m_C(A,B) \cong \Alg(C \triangleright A, B)$
	corresponds to the the algebra morphism
	$
		\circ_\m\colon \Alg(D \triangleright B, T)\times \Alg(C \triangleright A, B) \to \Alg(D \triangleright (C \triangleright A), T),
		(\tilde{\psi}, \tilde{f}) \mapsto \tilde{\psi} \circ (\id_D \triangleright \tilde{f}).
	$
	Using Lemma~\ref{lem:tensorcommutes}, we can state that $\Phi$ respecting composition of measurings is equivalent to verifying
	\[\begin{tikzcd}[row sep = tiny, column sep = small]
		{\Alg^F(D \triangleright B, T) \times \Alg^F(C \triangleright A, B)} && {\Alg^F(D \triangleright (C \triangleright A), T)} \\
		\\
		{\Alg^G(\mu_!(D \triangleright B), \mu_!(T)) \times \Alg^G(\mu_!(C \triangleright A), \mu_!(B))} && {\Alg^G(\mu_!(D \triangleright (C \triangleright A)), \mu_!(T))} \\
		&& {\Alg^G(\mu_*(D) \triangleright \mu_!(C \triangleright A), \mu_!(T))} \\
		{\Alg^G(\mu_*(D) \triangleright \mu_!(B), \mu_!(T)) \times \Alg^G(\mu_*(C) \triangleright \mu_!(A), \mu_!(B))} && {\Alg^G(\mu_*(D) \triangleright (\mu_*(C) \triangleright \mu_!(A)), \mu_!(T))}
		\arrow["{\circ_\m^F}", from=1-1, to=1-3]
		\arrow["{\mu_! \times \mu_!}"', from=1-1, to=3-1]
		\arrow["{\mu_!}", from=1-3, to=3-3]
		\arrow["\cong"', from=3-1, to=5-1]
		\arrow["\cong", from=3-3, to=4-3]
		\arrow["\cong", from=4-3, to=5-3]
		\arrow["{\circ_\m^G}", from=5-1, to=5-3]
	\end{tikzcd}\]
	commutes.
	Given a pair $(\tilde{\psi}, \tilde{f})\in \Alg^F(D \triangleright B, T) \times \Alg^F(C \triangleright A, B)$, we need to verify
	$
	\mu_!(\tilde{\psi} \circ (\id_D \triangleright \tilde{f})) = \mu_!(\tilde{\psi}) \circ \mu_!(\id_D \triangleright \tilde{f})
	$
	corresponds to
	$
	\mu_!(\tilde{\psi}) \circ (\id_{\mu_*(D)} \triangleright \mu_!(\tilde{f}))
	$
	under the isomorphism $\Alg^G(\mu_!(D \triangleright (C \triangleright A)), \mu_!(T)) \cong \Alg^G(\mu_*(D) \triangleright \mu_!(C \triangleright A), \mu_!(T))$.
	This is indeed the case by naturality of the isomorphism $\mu_!(D \triangleright (C \triangleright A)) \cong \mu_*(D) \triangleright \mu_!(C \triangleright A)$.
\end{proof}

\noindent\textbf{Lemma 3.24}
\textit{	
	The natural transformation defined in Definition~\ref{def:familyoffunctions2} $\Phi$ respects composition.}

\begin{proof}
	We know $\Phi(f) = f \circ (\epsilon_C \otimes \id_A)$, so we explicitly check composition of measurings is respected.
	Given $G$-meausurings $f \colon C \otimes A \to B$ and $\psi \colon D \otimes B \to T$, it is necessary to check the diagram
	\[\begin{tikzcd}[column sep=8.5em]
		{\mu_{\invexcl}(D) \otimes \mu_{\invexcl}(C) \otimes A} & {\mu_{\invexcl}(D) \otimes B} & T \\
		{\mu_{\invexcl}(D \otimes C) \otimes A}
		\arrow["{\id_{\mu_{\invexcl}(D)} \otimes (f \circ (\epsilon_C \otimes \id_A))}", from=1-1, to=1-2]
		\arrow["{\nabla^{\mu_{\invexcl}}}"', from=1-1, to=2-1]
		\arrow["{\psi \circ (\epsilon_D \otimes \id_B)}", from=1-2, to=1-3]
		\arrow["{(\psi \circ (\id_D\otimes f)) \circ (\epsilon_{D\otimes C} \otimes \id_A)}"', from=2-1, to=1-3]
	\end{tikzcd}\]
	commutes.
	This is the case since $\epsilon \colon \mu_* \circ \mu_{\invexcl} \to \id$ is a monoidal natural transformation, hence as morphisms in $\C$ we have
	$
	\epsilon_D\otimes \epsilon_C = \epsilon_{C \otimes D} \circ \nabla^{\mu_{\invexcl}}.
	$
\end{proof}

\end{document}